\documentclass{gtart_a}     
\pdfoutput=1

\usepackage{pinlabel}

\title[Canonical triangulations]{On canonical triangulations of 
once-punctured\\torus bundles and two-bridge link complements}

\author[Fran\c{c}ois Gu\'eritaud and David Futer]{Fran\c{c}ois Gu\'eritaud\newline Appendix by David Futer}
\givenname{Fran\c{c}ois}
\surname{Gu\'eritaud}
\address{DMA, \'Ecole normale sup\'erieure, CNRS\\\newline
45 rue d'Ulm\\75005 Paris\\France\\\newline\\\newline
{\rm[D Futer]}\qua Mathematics Department\\
Michigan State University\\\newline
East Lansing, MI 48824\\USA}
\email{gueritau@dma.ens.fr}
\urladdr{}

\givenname{David}
\surname{Futer}
\email{dfuter@math.msu.edu}
\urladdr{}

\volumenumber{10}
\issuenumber{}
\publicationyear{2006}
\papernumber{30}
\lognumber{0677}
\startpage{1239}
\endpage{1284}

\doi{}
\MR{}
\Zbl{}

\arxivreference{math.GT/0406242}


\received{10 November 2005}
\revised{29 July 2006}
\accepted{23 July 2006}
\published{16 September 2006}
\publishedonline{16 September 2006}
\proposed{Walter Neumann}
\seconded{Jean-Pierre Otal, Joan Birman}
\corresponding{}
\editor{Walter Neumann}
\version{}

\AtBeginDocument{\def\leftmark{Fran\c{c}ois Gu\'eritaud}}
\AtBeginDocument{\let\bar\wbar\let\tilde\wtilde}
\makeop{PSL}
\makeop{SL}
\newcommand{\smatris}[4]{\bigl( \begin{smallmatrix}
#1 & #2 \\ #3 & #4 \end{smallmatrix} \bigr)}

\makeatletter
\def\cnewtheorem#1[#2]#3{\newtheorem{#1}{#3}[section]
\expandafter\let\csname c@#1\endcsname\c@theorem}

\newtheorem{theorem}{Theorem}[section]

\cnewtheorem{lemma}[theorem]{Lemma}
\cnewtheorem{sublemma}[theorem]{Sublemma}
\cnewtheorem{proposition}[theorem]{Proposition}
\cnewtheorem{corollary}[theorem]{Corollary}

\theoremstyle{remark}
\cnewtheorem{trivial}[theorem]{Easy Fact}
\cnewtheorem{definition}[theorem]{Definition}

\makeatother  

\newcommand{\matris}[4]{\left ( \!\! \begin{array}{cc}
#1 & #2 \\ #3 & #4 \end{array} \!\! \right )}

\newcommand{\var}{\varepsilon}

\newcommand{\llabel}[1]{\label{#1}}

\renewcommand{\setminus}{{\smallsetminus}}
\newcommand{\cross}{{\times}} 
\newcommand{\bdy}{{\partial}} 
\newcommand{\vol}{{\mathcal{V}}}
\newcommand{\tw}{{\mathrm{tw}}}

\begin{document}

\begin{asciiabstract}
We prove the hyperbolization theorem for punctured torus bundles and
two-bridge link complements by decomposing them into ideal tetrahedra
which are then given hyperbolic structures, following Rivin's volume
maximization principle.
\end{asciiabstract}

\begin{abstract}
We prove the hyperbolization theorem for punctured torus bundles and
two-bridge link complements by decomposing them into ideal tetrahedra
which are then given hyperbolic structures, following Rivin's volume
maximization principle.
\end{abstract}

\maketitle
\cl{\small\it\`A la m\'emoire de Pierre Philipps}

\section{Introduction}

Let $T:=({\mathbb R}^2 \setminus {\mathbb Z}^2)/{\mathbb Z}^2$ be the once-punctured torus endowed with its differential structure and an orientation. The group $G$ of
isotopy classes of orientation-preserving diffeomorphisms $\varphi \co  T \rightarrow T$ (or the mapping class group of $T$) is identified as $G\simeq \SL_2({\mathbb Z})$, so each
such map $\varphi$ has well-defined eigenvalues in ${\mathbb C}$. For $[\varphi]$ in $G$, define the \emph{punctured torus bundle\/} $$V_{\varphi}:=T\times
[0,1]/\sim$$ where $\sim$ identifies $(x,0)$ with $(\varphi(x),1)$ for all $x$ in $T$. Then $V_{\varphi}$ is a differentiable oriented $3$--manifold, well-defined up to
diffeomorphism. Thurston's Hyperbolization Theorem \cite{otal} implies the following theorem as a very special case.

\begin{theorem} \llabel{main} If $\varphi$ has two distinct real eigenvalues, the punctured torus bundle $V_{\varphi}$ admits a finite-volume, complete hyperbolic
metric. \end{theorem}

The aim of this paper is to prove \fullref{main} by elementary and, to some extent, constructive arguments. The strategy is to exhibit a canonical, geodesic triangulation ${\mathcal H}$ of $V_{\varphi}$ into ideal tetrahedra (hyperbolic tetrahedra whose vertices are at infinity).

Combinatorially, ${\mathcal H}$ (sometimes called the Floyd--Hatcher or monodromy triangulation) is found by expressing a certain conjugate of $\pm \varphi$ as a product
of positive transvection matrices. Once such combinatorial data for a triangulation is given, the problem of making it hyperbolic lends itself to two approaches. One is
complex, explicit and ``local'': cross-ratio computations, particular hyperbolic isometries, etc (see eg\ Neumann and Zagier \cite{zagier}). The other approach, first described by 
Rivin \cite{rivin}, de Verdi\`ere
\cite{colin} and Casson, is real-analytic and ``global'': in order to make the structure complete, one kills its monodromy by maximizing the total hyperbolic volume (but
combinatorial obstructions might arise). In the case of $V_{\varphi}$, the combinatorial structure of ${\mathcal H}$ is sufficiently well-understood to allow a nice
interplay between the two approaches, yielding useful ``medium-range'' results (\fullref{ageometricallemma}). The philosophy of such results is that if the
structure with maximal volume is noncomplete, it should still be complete at ``most'' places, enabling us to make geometric statements.

Akiyoshi \cite{akiyoshi}, combining the methods of Akiyoshi, Sakuma, Wada and Yamashita \cite{aswy-announce} and Minsky \cite{minsky}, proved that the ``combinatorially canonical'' triangulation
${\mathcal H}$ must also be ``geometrically canonical,'' ie topologically dual to the Ford--Voronoi domain of $V_{\varphi}$. Lackenby \cite{lackenby}, assuming the existence of the hyperbolic metric, derived the same result by a normal surface argument. In \cite{quasifuchsien}, we apply
the methods of the present paper to extend the Akiyoshi--Lackenby theorem to quasifuchsian groups (where pleating laminations of the convex core replace the attractive and repulsive laminations of the monodromy $\varphi$). 

Knowing the space of angle structures also allows for easy volume estimates, some of which are worked out in the Appendix: these estimates make the constants of Brock \cite{brock} explicit (and sharp) for the class of manifolds under consideration. Although the main results of the present paper are known, our ambition is to demonstrate that hyperbolic geometry and combinatorics (of the curve complex, say) can interact more intimately than at the level of coarse geometry, a phenomenon which seems to extend beyond punctured torus groups and begs to be further explored. Other references closely related to this subject include the beautiful article of Akiyoshi, Sakuma, Wada and Yamashita \cite{aswy}, which builds on the work of J{\o}rgensen and partly motivated our work \cite{quasifuchsien}, and the examples compiled in 
Alestalo and Helling \cite{sfb343-1}, Helling \cite{sfb343-2} and Koch \cite{koch}.

The converse of \fullref{main} is true. If the trace $\tau$ of the monodromy map $\varphi$ is in $\{-1,0,1\}$, then $[\varphi]$ has finite order and $V_{\varphi}$ is Seifert fibered. If $\tau=\pm 2$, then $\varphi$ preserves a nontrivial simple closed curve $\gamma$ (parallel to a rational eigenvector) in the punctured torus, and $\gamma$ defines an incompressible torus or Klein bottle in $V_{\varphi}$. In any case we get a topological obstruction to the existence of the hyperbolic metric.

An attempt to be self-contained will be made in proving \fullref{main}. The proof will deal primarily with the case where the eigenvalues of $\varphi$ are positive. The other case is only a minor variant (in particular, $V_{-\varphi}$ can be obtained by ungluing the metric tetrahedra of $V_{\varphi}$ and regluing them in a slightly different way).

\fullref{prelim} is standard and recalls the classification of conjugacy classes in $\SL_2({\mathbb Z})$ in order to define the triangulation ${\mathcal H}$. The latter is studied in greater detail in \fullref{topotorinf} and \fullref{combotorinf}. Positive angles for ${\mathcal H}$ (a ``linear hyperbolic structure'') are provided in \fullref{positiveangles}. \fullref{hyperbolicvolume} explains the role played by hyperbolic volume maximization, allowing us to deal with the easy
cases in \fullref{somemore}. \fullref{ageometricallemma} presents the essential geometric lemma for the final attack, carried out in \fullref{hallali}. \fullref{numericalexample} is devoted to a numerical example. In \fullref{spheres}, we quickly recall the connection between once-punctured tori and $4$--punctured spheres. In the Appendix, David Futer builds on that connection to prove an analogue of \fullref{main} for the complements of two-bridge links and computes geometric estimates. 

The symbol ``$=$'' is preceded (resp.\ followed) by a colon ($:=$, resp.\ $=:$) when the equality serves as a definition for its left (resp.\ right) member.

I would like to thank Francis Bonahon and Fr\'ed\'eric Paulin for numerous insights and the great improvements this paper owes to them. Exciting discussions with David Futer and with Makoto Sakuma also gave invaluable input. This paper reached its pre-final form during a stay at the Institut Bernoulli (EPFL, Lausanne) for whose kind hospitality I express my deep gratitude. Finally, thanks are due to the referee for many helpful comments. This work was partially supported by NSF grant DMS-0103511.

\section[Conjugacy in SL(2,Z) and the Farey tesselation]{Conjugacy in $\SL_2({\mathbb Z})$ and the Farey tesselation} \llabel{prelim}

\begin{proposition} Let $\varphi$ be an element of $\SL_2({\mathbb Z})$ with two distinct eigenvalues in ${\mathbb R}_+^*$. Then the conjugacy class of $\varphi$ in
$\SL_2({\mathbb Z})$ contains an element of the form $$A\varphi A^{-1}= \matris{1}{a_1}{0}{1} \! \matris{1}{0}{b_1}{1} ~~ \matris{1}{a_2}{0}{1} \! \matris{1}{0}{b_2}{1}
~ \dots ~ \matris{1}{a_n}{0}{1} \! \matris{1}{0}{b_n}{1}$$ where $n>0$ and the $a_i$ and $b_i$ are positive integers. Moreover, the right hand side is unique up to
cyclic permutation of the factors $\smatris{1}{a_i}{0}{1} \! \smatris{1}{0}{b_i}{1} $. Conversely, any nonempty product of such factors is an element of $\SL_2({\mathbb
Z})$ with two distinct eigenvalues in ${\mathbb R}_+^*$. \llabel{classif} 
\end{proposition}

We sketch a proof of this popular fact, mainly in order to introduce the \emph{cyclic word\/} $\Omega$ associated to $\varphi$. The converse implication is easy (just
check that the trace is larger than $2$), so we focus on the direct statement.

Consider the upper half-plane model of the hyperbolic plane ${\mathbb H}^2$, endowed with the Farey tessellation $F$ (the ideal triangle $01\infty$ iteratedly
reflected in its sides). We identify $\PSL_2(\mathbb{R})$ with the group of isometries of $\mathbb{H}^2$ \emph{via\/} the isomorphism $\Psi$ defined
by $$\Psi\matris{a}{b}{c}{d}\co z \mapsto \frac{dz+c}{bz+a}.$$ (Under this slightly unusual convention, the slopes of the eigenvectors of $M$ are the fixed points of
$\Psi(M)$, rather than their inverses as would normally be the case.) It is known that the group of orientation-preserving isometries of ${\mathbb H}^2$ preserving $F$
is thus identified with $\PSL_2({\mathbb Z})$.

If $D$ is the oriented hyperbolic line running from the repulsive fixed point of $\pm\varphi$ to the attractive one, then $D$ crosses infinitely many Farey triangles
$(\dots t_{-1},t_0,t_1,t_2,\dots)$ of $F$. We can formally write down a bi-infinite word $$\Omega=\dots LRRRLLR\dots$$ where the $k$--th letter is $R$ (resp.\ $L$) if
$D$ exits $t_k$ to the Right (resp.\ Left) of where it enters. We will also say that $D$ \emph{makes a Right\/} (resp.\ \emph{Left\/}) at $t_k$. The word $\Omega$ contains
at least one $R$ and one $L$, because the ends of $D$ are distinct. The image of $t_0$ under $\varphi$ is a certain $t_m$ ($m>1$), and $\Omega$ is periodic of period
$m$.

Next, define the standard transvection matrices $$R:=\matris{1}{1}{0}{1} \quad \text{and} \quad L:=\matris{1}{0}{1}{1}.$$ These are parabolic transformations of ${\mathbb H}^2$ whose
respective fixed points are $0$ and $\infty$. Let $M$ be any subword of $\Omega$ of length $m$: we see $M$ as a product of standard transvection matrices, and
therefore as an element of $\SL_2(\mathbb{Z})$. By studying the actions of $R$ and $L$ on $F$, it is then easy to see that $\varphi$ and $M$ are conjugates in
$\PSL_2(\mathbb{Z})$, and therefore in $\SL_2(\mathbb{Z})$ since both have positive trace. This proves the existence statement for the $(a_i,b_i)$.

Uniqueness is checked as follows: on one hand, if $\varphi$ and $\varphi'$ are conjugates, a certain element of $\PSL_2({\mathbb Z})$ (preserving $F$) takes the axis of $\varphi$ to the axis of $\varphi'$, so they define the same word $\Omega$ up to translation. On the other hand, looking at the actions of $R$ and $L$ on ${\mathbb H}^2$, one
sees that a product of standard transvection matrices (as in the statement of \fullref{classif}) will define the word $\Omega=R^{a_1}L^{b_1}\dots
R^{a_n}L^{b_n}$, concatenated infinitely many times with itself.

In the language of \fullref{classif}, the sequence $(a_1,b_1,\dots,a_n,b_n)$ can be shown to be (a positive power of) the period of the continued fraction
expansion of the slope of the expanding eigenvector of $\varphi$. The word $\Omega$ will be seen either as infinite periodic, or as finite cyclic, depending on the
context.

\section{The canonical triangulation} \llabel{topotorinf}

\subsection{Diagonal exchanges} \llabel{leucocephale} There is another well-known interpretation of the Fa\-rey tessellation $F$ of the hyperbolic plane ${\mathbb H}^2$.
Under the canonical identification $H_1(T,{\mathbb Z})\simeq {\mathbb Z}^2$, where $T$ is the punctured torus defined in the Introduction, each rational number in the
boundary $\widehat{\mathbb R}={\mathbb R}\cup \{\infty\}$ of ${\mathbb H}^2$ can be seen as a \emph{slope\/}, ie a proper isotopy class of properly embedded lines in $T$,
going from the puncture to itself. The action of $\SL_2({\mathbb Z})$ on $\widehat{\mathbb Q}$ coincides with the action of the mapping class group $G$ of $T$ on rational
slopes. The edges of the Farey tessellation $F$ connect exactly the pairs of rational numbers whose corresponding slopes, or curves in $T$, can be homotoped off each other (away from the
puncture). The faces of $F$, having three edges, correspond exactly to the isotopy classes of \emph{ideal triangulations\/} of $T$: any such triangulation has one vertex (the
puncture), three edges, and two triangles (which meet along each edge). As one crosses from a face of $F$ to one of its neighbors, exactly one vertex gets replaced, which in
the corresponding triangulations of $T$ means that exactly one edge is changed. Inspection shows that the triangulation must be undergoing a \emph{diagonal exchange\/}: erase
one edge $e$, thus liberating a quadrilateral space $Q$ of which $e$ was a diagonal, then insert back the other diagonal.

\subsection{Tetrahedra} As before, let $\varphi$ be an element of $\SL_2({\mathbb Z})$ with two distinct eigenvalues in ${\mathbb R}_+^*$. In the proof of \fullref{classif}, we introduced the triangles $t_0,t_1,\dots$ crossed by the axis $D$ of $\varphi$. In view of the above, this yields a nonbacktracking path of diagonal
exchanges between some triangulation (associated to $t_0$) and its pushforward by $\varphi$ (associated to $t_m$).

More precisely, when the oriented line $D$ crosses an edge $e$ of the Farey tessellation, $e$ comes with a transverse orientation. So we can define the \emph{top\/}
(resp.\ \emph{bottom\/}) triangulation $\tau_+(e)$ (resp.\ $\tau_-(e)$) of the punctured torus $T$ as being the one associated with the Farey triangle crossed just after
(resp.\ before) the edge $e$. A diagonal exchange separates the triangulations $\tau_-(e)$ and $\tau_+(e)$. An ideal  tetrahedron is by definition a space diffeomorphic
to an ideal hyperbolic tetrahedron (topologically it is a compact tetrahedron with its vertices removed). We can immerse such an ideal tetrahedron $\Delta(e)$ in
$T\times {\mathbb R}$: the boundary of the immersed $\Delta(e)$ is made up of two \emph{pleated surfaces\/} (top and bottom) homotopic to $T$ and triangulated according
to $\tau_+(e)$ and $\tau_-(e)$ respectively (\fullref{onetetrahedron}). The immersion is an embedding on the interior of $\Delta(e)$, but two pairs of opposite
edges undergo identifications.
\begin{figure}[ht!] \centering \includegraphics{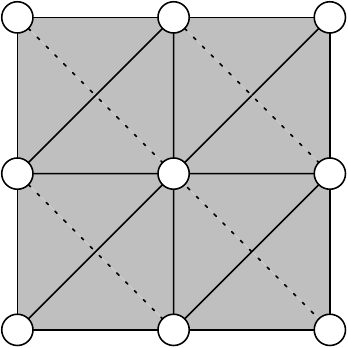} \caption{Four copies of $\Delta(e)$ in the cover $({\mathbb R}^2\setminus {\mathbb Z}^2)\times{\mathbb R}$ of
$T\times {\mathbb R}$} \llabel{onetetrahedron} \end{figure}

Next, if $D$ crosses the Farey edges $e_i,e_{i+1},\dots$, we can glue the top of the tetrahedron $\Delta_i:=\Delta(e_i)$ onto the bottom of $\Delta_{i+1}$ in
$T\times{\mathbb R}$, because $\tau_+(e_i)=\tau_-(e_{i+1})$. We thus get a bi-infinite stack of tetrahedra $(\Delta_i)_{i\in{\mathbb Z}}$. For any nonnegative $N$ the space
$U_N:=\bigcup_{i=-N}^N \Delta_i$ is a strong deformation retraction of $T\times {\mathbb R}$. For $N$ large enough, $U_N$ is homeomorphic to  $T\times [0,1]$:
the only way this can fail is if all the $\Delta_i$ for $-N\leq i \leq N$ have a common edge; but any edge of any tetrahedron $\Delta_j$ is shared by only finitely
many other (consecutive) $\Delta_i$'s, because for any Farey vertex $v$, only finitely many of the Farey edges $e_i$ bound triangles with $v$ as a vertex (and these
$e_i$ are consecutive). Therefore, the space $U=\bigcup_{i\in{\mathbb Z}}\Delta_i$ is homeomorphic to $T\times {\mathbb R}$. If  $m$ is the period of the word
$\Omega$, there is an orientation-preserving homeomorphism $\Phi$ of $U$, acting like $[\varphi]$ on the $T$--factor, that sends $\Delta_i$ to $\Delta_{i+m}$ for all
$i$. The quotient $U/\Phi$ is a manifold (homeomorphic to) $V_{\varphi}$, naturally triangulated into $m$ ideal  tetrahedra.

\fullref{whitehead} also shows a way to interpret the standard transvection matrices $R$ and $L$ directly as adjunctions of new tetrahedra (by performing diagonal exchanges on the top faces). Similarly, to topologically triangulate a general pseudo-Anosov surface bundle, we can always go by diagonal exchanges from some triangulation (of the surface) to its pushforward by the monodromy map, an idea usually credited to Casson.
\begin{figure}[ht!] \centering \includegraphics{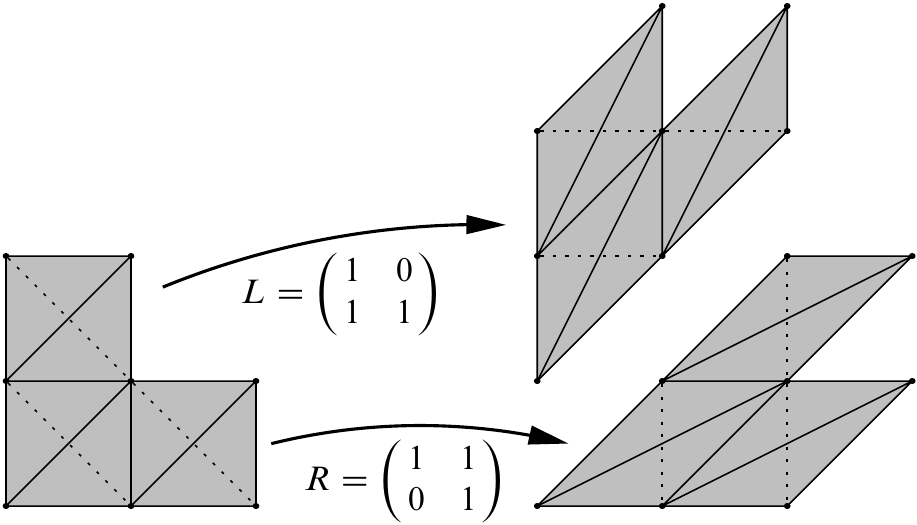} \caption{} \llabel{whitehead} \end{figure}

\section{Combinatorics of the torus at infinity} \llabel{combotorinf}

The manifold $V_{\varphi}$ is naturally homeomorphic to the interior of a compact manifold with boundary, denoted by $\overline{V}_{\varphi}$ and defined as a bundle over
the circle with fiber $T-\delta$, where $\delta$ is a regular neighborhood of the puncture.

The \emph{torus at infinity\/} of the manifold $V_{\varphi}$ is the boundary of $\overline{V}_{\varphi}$, namely a topological torus. The links of the vertices of the
tetrahedra $\Delta_i$ provide a tessellation ${\mathcal A}$ of the torus at infinity into topological triangles. In this section we investigate the combinatorics of
${\mathcal A}$.

Each vertex of ${\mathcal A}$ corresponds to an edge of $V_{\varphi}$ shared by a few consecutive tetrahedra $\Delta_i$. This edge in turn corresponds to a Farey
vertex shared by a few consecutive Farey triangles. The union of all the Farey triangles adjacent to a given vertex $v$ forms a \emph{fan\/}. If $v$ arises as a vertex
of triangles visited by the oriented line $D$, one of the following two things must happen right after $D$ enters the fan: either $D$ makes a Right, then a number of Lefts
(possibly 0), then a Right and leaves the fan; or the same is true, exchanging Right and Left.

Therefore, the vertices of ${\mathcal A}$ correspond exactly to the subwords of $\Omega$ of the form $RL^*R$ or $LR^*L$ (where $*\geq 0$). Each such subword actually
corresponds to two vertices of ${\mathcal A}$, because the edges of the tetrahedra $\Delta_i$ have two ends.
\begin{figure}[ht!] \centering \includegraphics{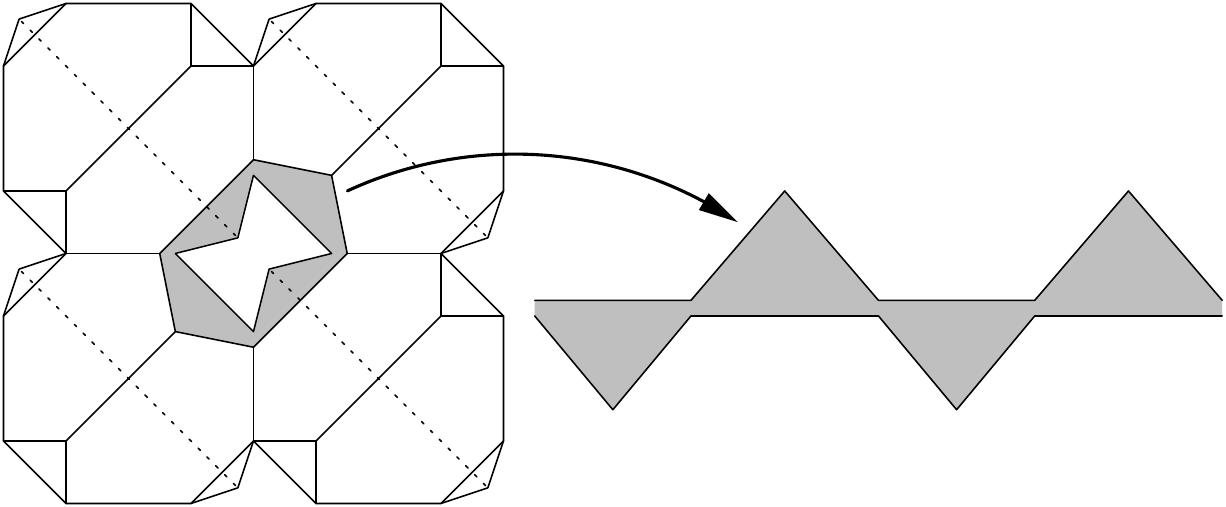} \caption{The link of the puncture} \llabel{fourtriangles} \end{figure}

Moreover, each tetrahedron $\Delta_i$, having four vertices, contributes exactly four triangles to ${\mathcal A}$. By looking at a vertex (puncture) of the cover
$({\mathbb R}^2 \setminus {\mathbb Z}^2)\times {\mathbb R}$ of $T\times {\mathbb R}$ with embedded $\Delta_i$, one checks (\fullref{fourtriangles}) that each of the
four triangles has exactly one vertex not shared with any of the other three: we call this vertex the \emph{apex\/} and the opposite edge the \emph{base\/}. The four bases
form a broken line of four segments which is a closed curve running around the puncture, and the four apices point alternatively up and down in the ${\mathbb
R}$--factor. Such chains of four triangles must be stacked on top of one another while respecting the previously described combinatorics of the vertices. The result is
shown in \fullref{pile}, where the underlying word $\Omega$ was chosen to be $\dots R^4L^4R^4L^4\dots$ (read from bottom to top). A few remarks are in order.
\begin{figure}[ht!] \centering \includegraphics{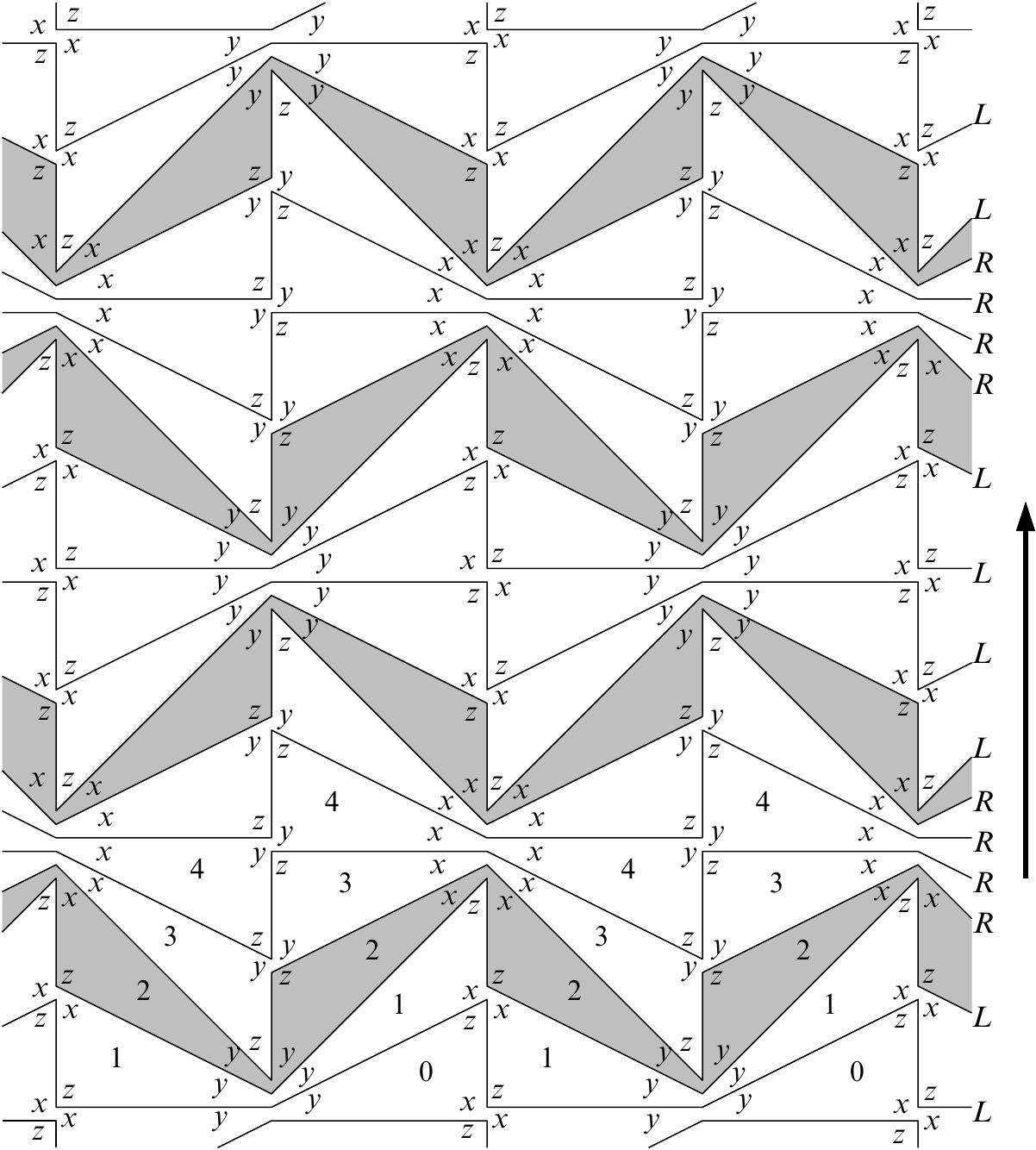} \caption{The tessellation $\mathcal{A}$} \llabel{pile} \end{figure}

\medskip \textbf{First remark}\qua We labeled by $x_i,y_i,z_i$ the angular sectors of the triangles corresponding to the tetrahedron $\Delta_i$ (in \fullref{pile},
the subscript $i$ is replaced by a number at the center of the triangle, omitted after the first few levels). Each angular sector corresponds to a (topological)
dihedral angle of $\Delta_i$. Opposite dihedral angles are equal in ideal hyperbolic tetrahedra: this is why three different labels per level are enough (instead of
six, the total number of edges in a tetrahedron).

\medskip \textbf{Second remark}\qua The design in \fullref{pile} of the vertices of the torus at infinity, represented as being ``opened up'', is intended to emphasize
the layered structure of the picture (each layer corresponds to one tetrahedron $\Delta_i$).

\medskip \textbf{Third remark}\qua Proving \fullref{main} by the method outlined in the Introduction amounts to realizing \fullref{pile} geometrically in the plane by Euclidean triangles, making same-layer
angles with identical labels equal (\fullref{rivinlemma} will make this statement more precise).

\medskip \textbf{Fourth remark}\qua The convention is that the pair of edges of $\Delta_i$ that does not get identified into one edge is labeled $z_i$: so $z_i$ is the label at the apex (in the sense of \fullref{fourtriangles}). Equivalently, if a tetrahedron is seen as a diagonal exchange, $z_i$ is the label common to the appearing and disappearing edges. The other edge pairs are labeled in such a way that in each triangle the letters $x,y,z$ appear \emph{clockwise\/} in that order. Therefore, if $e_i$ is a Farey edge and $p$ (resp.\ $q$) its right (resp.\ left) end for the transverse orientation, the dihedral angle of the tetrahedron $\Delta_i$ at the edge of slope $p$ (resp.\ $q$) is $x_i$ (resp.\ $y_i$).

\medskip \textbf{Fifth remark}\qua The tetrahedra $\Delta_i$ are naturally indexed in ${\mathbb Z}/m{\mathbb Z}$. The letters $R$ and $L$ live naturally on the pleated surfaces, or \emph{between\/} the tetrahedra $\Delta_i$. In \fullref{pile} and henceforth, the $i$--th layer is colored in grey if and only if the letters just before and just after $\Delta_i$ are different (here $i=2,6,\dots$). Such indices $i$ are called \emph{hinge\/} indices, because they are at the hinge between two nonempty subwords $R^p$ and $L^q$. Hinge tetrahedra (the associated $\Delta_i$) will play an important role later on.

\medskip \textbf{Sixth remark}\qua While the fundamental domain of \fullref{pile} is supposed to have a horizontal length of four triangles (see \fullref{fourtriangles}), we notice a horizontal period of length only two. This corresponds to the ``hyperelliptic'' involution of the once-punctured fiber torus (rotation of $180^{\circ}$ around the puncture, central in $\SL_2({\mathbb Z})$ and therefore well-defined on $V_{\varphi}$). This will simplify many of our computations.

\medskip \textbf{Seventh remark}\qua The valence of a vertex $s$ corresponding to a subword $RL^nR$ or $LR^nL$ of $\Omega$ (where $n\geq 0$) is $2n+4$. This is because exactly $n+2$ Farey triangles are adjacent to the corresponding rational $v_s$; each such Farey triangle defines a triangulated surface (with an edge of slope $v_s$), and each such surface contributes exactly two segments issuing from $s$ in $\mathcal{A}$.

\section{Finding positive angles} \llabel{positiveangles}

The tetrahedra $\Delta_i$ and $\Delta_{i+1}$ have two common triangular faces whose union forms a \emph{pleated punctured torus\/} $\Sigma$ properly isotopic to $T\times\{*\}$ in $T\times \mathbb{R}$. Moreover, $\Sigma$ receives a transverse ``upward'' orientation from the $\mathbb{R}$--factor. Suppose all tetrahedra $\Delta_i$ are endowed with dihedral angles. Let $e$ be an edge of $\Sigma$: if the sum of all dihedral angles at $e$ \emph{below\/} $\Sigma$ is $\pi+\alpha$, we call $\alpha$ the \emph{pleating angle\/} of
$\Sigma$ at $e$.

In this section we find positive dihedral angles for the ideal tetrahedra $\Delta_i$. More precisely, we describe the convex space $\Pi$ of positive angles $x_i,
y_i, z_i$ for the $\Delta_i$ such that:
\begin{equation} \left \{ \begin{array}{rl}
\text{i ---} & \text{For each $i$ in $\mathbb{Z}/m\mathbb{Z}$ one has $x_i+y_i+z_i=\pi$;} \\
\text{ii ---} & \text{The dihedral angles around any edge add up to $2\pi$;} \\ 
\text{iii ---} & \text{For each $i$ in $\mathbb{Z}/m\mathbb{Z}$, the three pleating angles of the pleated}\\ & \text{punctured torus between $\Delta_i$ and $\Delta_{i+1}$
add up to $0$.} \end{array} \right . \llabel{coherence} \end{equation}
Condition \eqref{coherence}-ii is necessary, though not sufficient, for a hyperbolic structure at the edges; Condition \eqref{coherence}-iii is necessary, though not sufficient,
to make the loop around the puncture of $T$ a parabolic isometry of $\mathbb{H}^3$ (see \fullref{computhol}). This ``cusp condition'' \eqref{coherence}-iii restricts the space of angle structures, but will make it a little easier to describe.

Recall the line $D$ that runs from the repulsive fixed point $q^-$ to the attractive fixed point $q^+$ of $\varphi$ across the Farey triangulation. If the tetrahedron
$\Delta_i$, corresponding to the Farey edge $e_i$, realizes a diagonal exchange that erases an edge $\var'$ and replaces it with $\var$, we denote by $z_i$ the interior
dihedral angle of $\Delta_i$ at $\var$ and $\var'$ by the fourth remark following \fullref{pile}. Observe that the slope of $\var$ (resp.\ $\var'$) is the rational located opposite
$e_i$ in the Farey diagram, on the side of $q^+$ (resp.\ $q^-$). We define the \emph{half pleating angle\/} $w_i$ by $\pi-2w_i=z_i$.

Thus, if \eqref{coherence}-ii and \eqref{coherence}-iii are to be satisfied, the pleating angles of the pleated punctured torus $\Sigma$ living between $\Delta_i$ and
$\Delta_{i+1}$ must be \begin{equation} \llabel{pleatingangles} -2w_i, \quad 2w_{i+1} \quad \text{and} \quad 2w_i-2w_{i+1}.\end{equation} (Observe the signs: by our definition,
angles pointing upwards, like the ``new'' edge of $\Delta_i$, are negative pleatings, while angles pointing downwards, like the ``old'' edge of $\Delta_{i+1}$, are
positive ones.) We write the numbers \eqref{pleatingangles} in the corners of the corresponding Farey triangle (top of \fullref{rl2}), distinguishing the cases
$L$ and $R$, and we repeat this for all indices $i$.
\begin{figure}[ht!] \centering \includegraphics{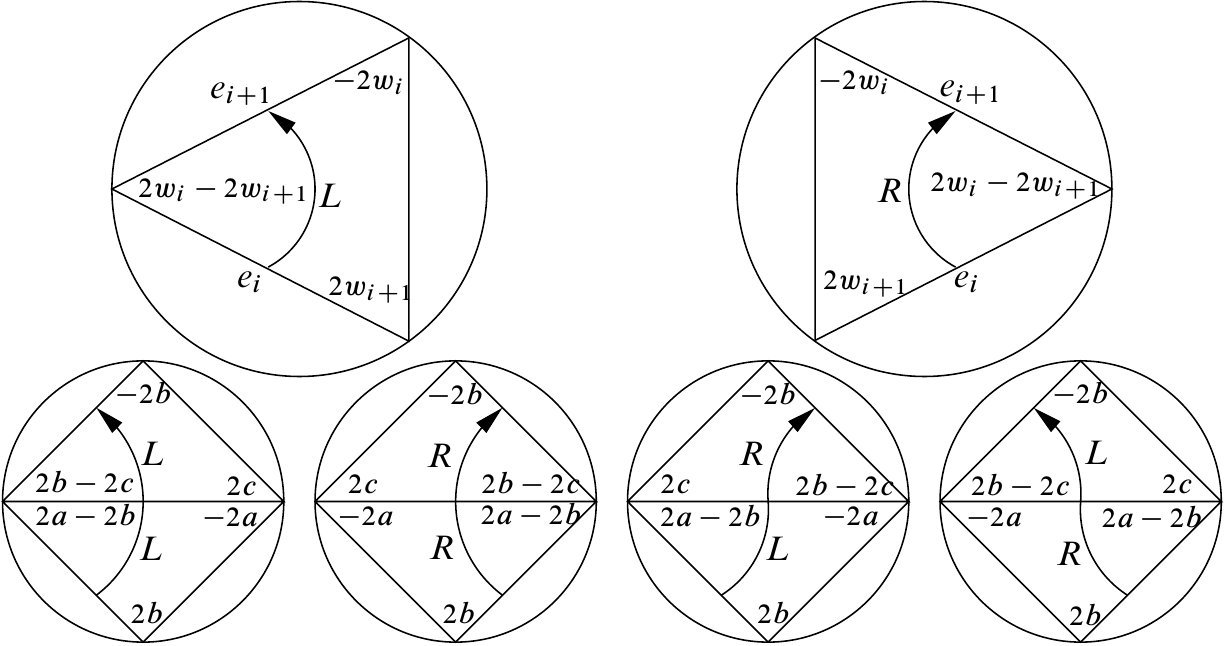} \caption{On the bottom row in each case, $e_i$ is the horizontal edge} \llabel{rl2} \end{figure}

In fact, the values of the $w_i$ will also determine the $x_i$ and $y_i$. To see this, we must write down the pleating angles of \emph{two\/} pleated surfaces, living
above \emph{and\/} below the tetrahedron $\Delta_i$. For notational convenience, write $$(w_{i-1},w_i,w_{i+1})=:(a,b,c).$$ By the fourth remark following \fullref{pile},
the quantity $2x_i$ (resp.\ $2y_i$) is the difference between the numbers written just above and just below the right (resp.\ left) end of $e_i$ in \fullref{rl2}
 (the factor $2$ comes from the fact that the two edges of $\Delta_i$ with angle $x_i$ (resp.\ $y_i$) are identified).  By computing differences between the pleating
angles given in \fullref{rl2}\,(bottom), we find the formulae in \fullref{interlettres} for $x_i, y_i, z_i$ (depending on the letters $\Omega_i^-$ and $\Omega_i^+$, each equal to
$R$ or $L$, living respectively just before and just after the index $i$).
\begin{table}[htb!] \llabel{interlettres} $$\begin{array}{c|c|c|c|c}
\Omega_i^-,\Omega_i^+&L~~~~~L&R~~~~~R&L~~~~~R&R~~~~~L\\ \hline
           x_i       &  a+c  &-a+2b-c& a+b-c &-a+b+c \\
           y_i       &-a+2b-c&  a+c  &-a+b+c & a+b-c \\
           z_i       & \pi-2b& \pi-2b&\pi-2b &\pi-2b \\\end{array}$$
\caption{}
\end{table}

Condition \eqref{coherence}-i can be checked immediately, while \eqref{coherence}-ii-iii are true by construction. From \fullref{interlettres}, the conditions for all
angles $x_i, y_i, z_i$ to be positive are that:
\begin{equation} \llabel{positivity} \left \{ \begin{array}{rl} \text{ i ---}&\text{For all $i$ one has $0<w_i<\pi/2$.} \\ \text{ ii ---}&\text{If $i$ is not a hinge
index (first two cases), $2w_i>w_{i+1}+w_{i-1}$.} \\ \text{ iii ---}&\text{If $i$ is a hinge index (last two cases), $|w_{i+1}-w_{i-1}|<w_i$.} \end{array} \right .
\end{equation}
We call \eqref{positivity}-i the \emph{range\/} condition, \eqref{positivity}-ii the \emph{concavity\/} condition, and \eqref{positivity}-iii the \emph{hinge\/} condition.
The space $P$ of sequences $(w_i)_{i\in\mathbb{Z}/m\mathbb{Z}}$ which satisfy \eqref{positivity}, homeomorphic to the solution space $\Pi$ of \eqref{coherence}, is
clearly an open, convex polyhedron of compact closure in $\mathbb{R}^m$. Moreover, $P$ is nonempty: to exhibit a sequence $(w_i)$ in $P$, set $w_j=\pi/3$ when $j$ is
a hinge index, and complete the gaps between consecutive hinge indices $j<k$ with strictly concave subsequences taking their values in $[\pi/3,\pi/2)$, eg\
$w_i=\pi/3+\frac{(i-j)(k-i)}{(k-j)^2}$ for $j\leq i \leq k$ (indices are of course seen in ${\mathbb Z}$ for the evaluation). \fullref{dubya} shows the typical
graph of a sequence $(w_i)$ that satisfies all conditions of \eqref{positivity}. Finally, note that the formulae of \fullref{interlettres} are still valid when $\Omega$
is reduced to $RL$ or $LR$ (the letters $a$ and $c$ are just two names for the same parameter then).
\begin{figure}[ht!] \centering \includegraphics{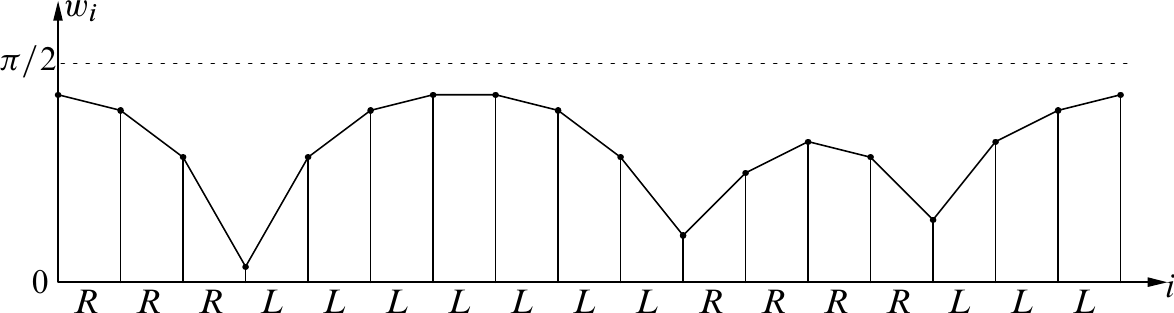} \caption{} \llabel{dubya} \end{figure}

\section{Hyperbolic volume} \llabel{hyperbolicvolume}

Our goal for the remainder of the paper is to find a point $(w_i)$ of $P$ where the tetrahedra fit together so as to yield a complete hyperbolic structure on
$V_{\varphi}$. This section is devoted to checking that this is the same as finding in $P$ a critical point of the total hyperbolic volume, an approach pioneered by
Rivin \cite{rivin}. A few facts concerning the volume of ideal tetrahedra will be needed.

\subsection{Volume of an ideal tetrahedron} The volume of a hyperbolic tetrahedron whose dihedral angles are $x,y,z>0$ is
\begin{equation}\llabel{untetraedre}\mathcal{V}(x,y,z)=-\int_0^x\log \,(2 \sin u)\,du - \int_0^y\log \,(2 \sin u)\,du - \int_0^z \log \,(2 \sin u)\,du\end{equation} (we refer to
Milnor \cite{milnorvolume} for a proof). Since $\int_0^{\pi}\log \,(2\sin u)\,du=0$, Equation \eqref{untetraedre} easily implies
the following proposition.

\begin{proposition} \llabel{tetvolume} The function $\mathcal{V}$ can be continuously extended by $0$ to all nonnegative triples $(x,y,z)$ such that $x+y+z=\pi$.   If
$\frac{d}{dt}(x_t,y_t,z_t)=(X,Y,Z)$ then \begin{equation} \exp \left ( \frac{-d}{dt} \mathcal{V}(x_t,y_t,z_t) \right ) =\sin^X x_t \sin^Y y_t
\sin^Z z_t. \llabel{tetvolumeformula} \end{equation} \end{proposition} \begin{proof} Straightforward. We will always apply this formula exactly in the form it is
stated, because it will usually make the right hand side simplest. \end{proof}

\subsection{Critical volume and trivial holonomy} 
\begin{lemma}[Rivin, Chan--Hodgson]\llabel{rivinlemma}  On the open affine polyhedron $P$ of sequences $(w_i)$ satisfying \eqref{positivity}, define the volume
functional $\vol$ as the sum of the hyperbolic volumes of tetrahedra $\Delta_i$ with dihedral angles $x_i, y_i, z_i$, as given by \fullref{interlettres}. Then
$(w_i)$ is a critical point for $\vol$ in $P$ if and only if the gluing of the tetrahedra $\Delta_i$ defines a complete finite-volume hyperbolic structure on the
punctured torus bundle $V_{\varphi}$. \end{lemma}
\begin{proof} This now standard fact holds for general ideal triangulations as well (see for example Chan \cite{chan-hodgson} or Rivin \cite{rivin}). The following proof, however, is
deliberately specific to the example at hand. This will enable us to introduce a few objects and relationships that will be useful in the sequel. Conversely, the main idea of the present proof (associate to each edge of $V_{\varphi}$ an explicit deformation in the space of angle structures) can be used to prove the general case of Rivin's theorem.

First we assume $(w_i)$ is critical. Let $B$ be the torus at infinity of $V_{\varphi}$ with the vertices of the tessellation ${\mathcal A}$ removed. Let $\sigma$ be the hyperelliptic involution of the fiber $T$, so $\sigma$ acts as a translation on $B$. Define $B'=B / \sigma$ and ${\mathcal A}' = {\mathcal A} / \sigma$. Let $t_0$ be a triangle of ${\mathcal
A}'$, $\epsilon_0$ an oriented edge of $t_0$ and $x_0$ an interior point of $t_0$. The group of orientation-preserving similarities of the Euclidean plane ${\mathbb C}$ is
${\mathbb C}^* \ltimes {\mathbb C}$.

\begin{definition} For a given $(w_i)$ in $P$, the \emph{holonomy function\/} is the representation $$\rho\co  \pi_1(B',x_0)  \rightarrow {\mathbb C}^* \ltimes {\mathbb C}$$
defined as follows. Given an element $\alpha$ of $\pi_1(B',x_0)$, view $\alpha$ as a cyclic sequence of triangles $t_0, t_1, \dots, t_s=t_0$ of ${\mathcal A}'$, such
that any two consecutive $t_i$'s share an edge. Then, draw an oriented copy $\tau_0$ of $t_0$ in the plane ${\mathbb C}$, with angles specified by $(w_i)$, by making
the image of the oriented edge $\epsilon_0$ coincide with $(0,1)$. Sharing an edge with $\tau_0$, draw a copy $\tau_1$ of $t_1$, also with angles specified by $(w_i)$. Then
draw a copy $\tau_2$ of $t_2$ adjacent to $\tau_1$, etc. By definition, $\rho(\alpha)$ is the orientation-preserving similarity mapping the copy of the oriented edge
$\epsilon_0$ in $\tau_0$ to the copy of $\epsilon_0$ in $\tau_s$. The \emph{reduced holonomy function\/} $\psi\co  \pi_1(B',x_0) \rightarrow {\mathbb C}^*$ is defined as the projection of
$\rho$ on the first factor. \end{definition}

It is a simple exercise to check that $\rho$ is well-defined, and is a representation (the concatenation rule being that $\alpha\beta$ denotes the path $\alpha$
followed by the path $\beta$). Note that $\psi$, having a commutative target, factors through a representation $\psi\co H_1(B',{\mathbb Z}) \rightarrow {\mathbb C}^*$.

\begin{sublemma} \llabel{computhol} Let $\alpha$ be an element of $H_1(B',{\mathbb Z})$ represented by a curve around a $4$--valent vertex of ${\mathcal A}'$, and let
$\beta$ be represented by a curve that follows a ``grey'' (hinge) level in ${\mathcal A}'$ (\fullref{holonomies}). If $(w_i)$ is critical for the volume functional
$\vol$, then $\psi(\alpha)=\psi(\beta)=1$. \end{sublemma}
\begin{figure}[ht!] \centering \includegraphics{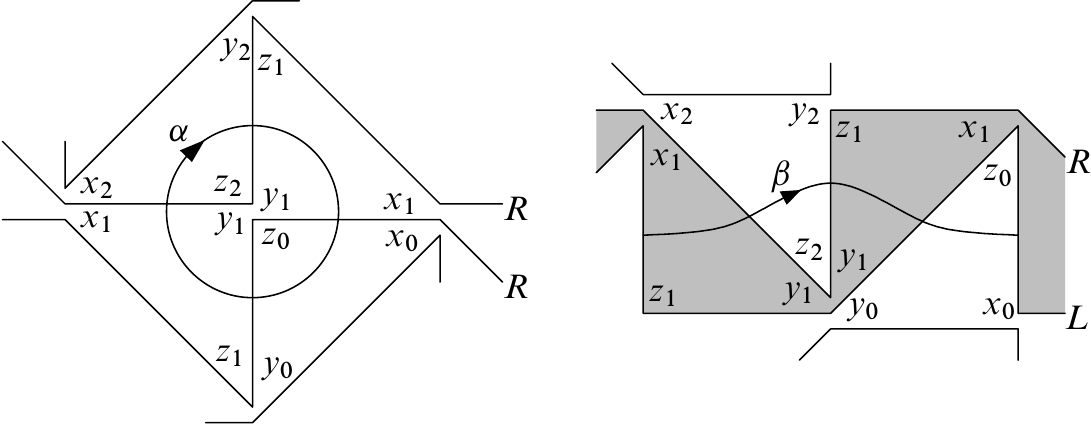} \caption{} \llabel{holonomies} \end{figure}
\begin{proof} We already know that $\psi(\alpha),\psi(\beta)$ belong to ${\mathbb R}_+^*$, since the  angle conditions \eqref{coherence} defining $P$ impose oriented
parallelism. It remains to prove that $|\psi(\alpha)|=|\psi(\beta)|=1$.

At a critical point, the partial derivatives of $\vol$ with respect to any of the $w_i$ must vanish. Between and near two identical letters, say $R$ and $R$,
according to \fullref{interlettres}, the angles are given by the table
$$\begin{array}{c|ccccc}
\Omega&      &R&       &R&       \\
   i  &   0  & &   1   & &   2   \\
 w_i  &   a  & &   b   & &   c   \\ \hline
 x_i  &\xi-b & &-a+2b-c& &-b+\xi'\\
 y_i  &\eta+b& &  a+c  & &b+\eta'\\
 z_i  &\pi-2a& & \pi-2b& & \pi-2c\end{array}$$
(the exact expression of $\xi,\xi',\eta,\eta'$ depends on the letters before and after $RR$, but only the $b$--contribution matters here). Using \fullref{tetvolume}, criticality of $\vol$ implies \begin{equation*} 1=\exp \frac{-\partial \vol}{\partial b}= \frac{\sin y_0\sin^2 x_1\sin y_2}{\sin x_0
\sin^2 z_1 \sin x_2}. \end{equation*} Using the fact that edge lengths in a triangle are proportional to angle sines, it follows that the edge lengths in \fullref{holonomies}\,(left) around the central vertex fit nicely together. So $|\psi(\alpha)|=1$. The case of a subword $LL$ is treated similarly, which takes care of
\emph{all\/} $4$--valent vertices of the tessellation ${\mathcal A}'$.

Near a hinge between two different letters, say $L$ followed by $R$, the angles are
$$\begin{array}{c|ccccc}
\Omega&      &L&      &R&       \\
   i  &  0   & &   1  & &   2   \\
 w_i  &  a   & &   b  & &   c   \\ \hline
 x_i  &\xi+b & &a+b-c & &-b+\xi'\\
 y_i  &\eta-b& &-a+b+c& &b+\eta'\\
 z_i  &\pi-2a& &\pi-2b& & \pi-2c\end{array}$$
This time, criticality implies \begin{equation*} 1=\exp \frac{-\partial \vol}{\partial b}=\frac{\sin x_0\sin y_1\sin x_1\sin y_2}{\sin y_0\sin^2 z_1\sin
x_2}. \end{equation*} As shown in \fullref{holonomies}\,(right) and by the same trigonometric argument, this means that the first and last edges crossed by $\beta$
have the same length. So $\psi(\beta)=1$. (If $\Omega$ is reduced to $LR$, the computation is formally the same, identifying indices $0$ and $2$.) The case of a
subword $RL$ is similar. \fullref{computhol} is proved. \end{proof}

Now let $\alpha$ be an element of $\pi_1(B',x_0)$ that is conjugated to a simple loop around a $4$--valent vertex of ${\mathcal A}'$. By \fullref{computhol} (and
an easy conjugation argument), $\rho(\alpha)$ is a translation. Moreover, $\rho(\alpha)$ fixes the vertex around which $\alpha$ revolves, so $\rho(\alpha)=1$, the
identity of the Euclidean plane.

\begin{sublemma} \llabel{trivialhol} Let $U$ be the quotient of the torus at infinity of $V_{\varphi}$ by the action of the hyperelliptic involution $\sigma$ of the fiber, so that $B'\subset U$. Suppose
$(w_i)$ is critical for $\vol$. Then the representation $\rho\co  \pi_1(B',x_0)  \rightarrow {\mathbb C}^* \ltimes {\mathbb C}$ descends to a representation
$\rho_U\co  \pi_1(U,x_0)  \rightarrow {\mathbb C}^*\ltimes {\mathbb C}$ whose first projection $\psi_U\co  \pi_1(U,x_0)  \rightarrow {\mathbb C}^*$ is trivial. \end{sublemma}
\begin{figure}[ht!] \centering \includegraphics{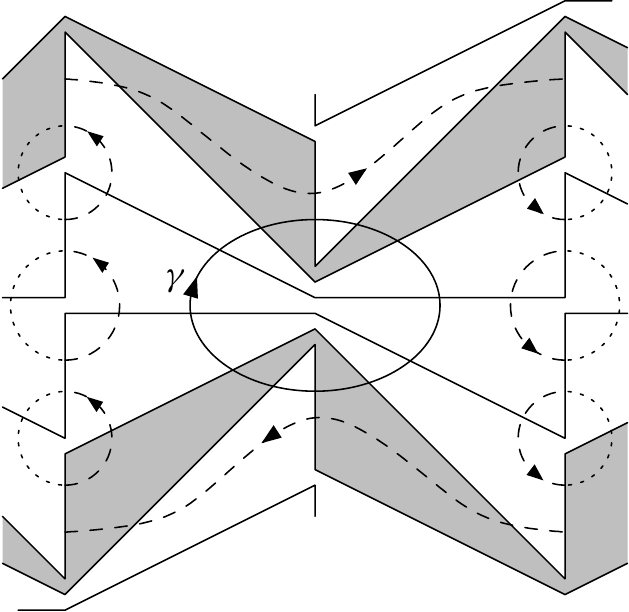} \caption{} \llabel{reconstruct}\end{figure}

\begin{proof} To see that $\rho_U$ is well-defined, we only need to check that, if $\gamma$ is (conjugated to) a loop around a vertex $v$ of ${\mathcal A'}$, then
$\rho(\gamma)=1$. If $v$ has valence $4$, it has already been done. If not, by the argument preceding \fullref{trivialhol}, it is sufficient to check that
$\psi([\gamma])=1$, where $[\gamma]$ denotes the homology class of $\gamma$. But in $H_1(B',{\mathbb Z})$, the element $[\gamma]$ is a sum of loops around $4$--valent
vertices and curves following ``grey'' levels (see \fullref{reconstruct}: the vertical edges of $B'$ on the two sides of the picture are identified, and the curves
crossing these edges undergo a ``split-and-paste'' process to yield $\gamma$). So by \fullref{computhol}, $\psi([\gamma])=1$; therefore $\rho_U$ is well-defined.
Moreover, if $\beta$ is a curve following a ``grey'' level, \fullref{computhol} tells that $\rho_U(\beta)=\rho(\beta)$ is a nonidentical Euclidean
\emph{translation\/}. The value of $\rho_U$ on another generator of $\pi_1(U)$ (which is abelian) must commute with $\rho_U(\beta)$, and therefore be a translation too.
So $\rho_U$ has its image contained in $\{1\}\ltimes{\mathbb C}$ and $\psi_U=1$, completing the proof. \end{proof}

By assigning length $1$, for example, to the reference edge $\epsilon_0$ of ${\mathcal A}'$, a critical point $(w_i)$ of the volume functional thus defines the lengths of all
other edges of ${\mathcal A}'$ in a coherent way. This yields a complete Euclidean metric $g$ on $U$. The universal cover $\widetilde{U}$ of $U$ thus embeds into
${\mathbb C}$ (the embedding, also called the developing map of the local Euclidean structure, can only be injective, because the $\tilde{g}$--geodesic joining two
distinct points of $\widetilde{U}$ is sent to a geodesic of $\mathbb{C}$); so there is an isometry $\widetilde{U} \simeq {\mathbb C}$. The metric $g$ lifts from $U$ to
the torus at infinity of $V_{\varphi}$ and its tessellation ${\mathcal A}$, producing a geometric realization of ${\mathcal A}$ and of \fullref{pile} in ${\mathbb C}$
(Euclidean plane tiling). Above each triangle of the universal cover of ${\mathcal A}$ now sits one ideal tetrahedron with vertex at infinity: the tetrahedron is the
hyperbolic convex hull of $\infty$ and the vertices of the triangle. Note that these tetrahedra fill ${\mathbb H}^3$ completely above a certain height.

To make sure that the pasted metric on the union $V=V_{\varphi}$ of all ideal tetrahedra is now complete, assume a geodesic $\gamma(t)$ in $V$ hits infinity at time $T<\infty$. If $K\subset V$ is compact, ie has a compact intersection with each tetrahedron $\Delta_i$, then $\gamma$ must eventually exit $K$ (if not, the $\gamma(T-1/n)$ accumulate at some point $p$ of some tetrahedron, but centered at $p$ there is a small embedded hyperbolic ball: absurd). So for $t$ sufficiently close
to $T$, there is a lift of $\gamma(t)$ arbitrarily high above the tessellation ${\mathcal A}$ (embedded in ${\mathbb C}$ in the upper half space model). But at
sufficiently great height, the tetrahedra above ${\mathcal A}$ fill ${\mathbb H}^3$ completely, so geodesics are defined for long times (eg\ times larger than $1$): a
contradiction. The first implication of \fullref{rivinlemma} is proved.

To prove the converse, it is enough to show that if the gluing of the tetrahedra yields a complete hyperbolic metric, then the gluing of their vertex links yields a
geometric realization of ${\mathcal A}$, ie of \fullref{pile} (checking $\partial\vol/\partial w_i=0$ then amounts to a rerun of the two computations in
\fullref{computhol}, distinguishing whether $i$ is a hinge index or not). But the latter is clear: given a complete hyperbolic metric, consider a triangulated
universal cover ${\mathbb H}^3\rightarrow V_{\varphi}$ and send (a lift of) the cusp to infinity in the upper half space model. It is a classic argument that deck
transformations of ${\mathbb H}^3$ which fix infinity must be parabolic, so the link of infinity has two translational periods and provides a Euclidean realization of
$\mathcal{A}$ (and of \fullref{pile}). \end{proof}

\subsection{Behavior of the volume functional} As a consequence of \fullref{rivinlemma}, to prove \fullref{main} we only need to find a critical point of the
volume functional $\vol$ in the open polyhedron $P$ of cyclic sequences $(w_i)$ satisfying the conditions \eqref{positivity}. A few more facts about the volume
of ideal hyperbolic tetrahedra will be needed.

By \fullref{tetvolume}, the volume functional $\mathcal{V}$ continuously extends to the (compact) closure $\bar{P}$ of the polyhedron $P$ (the space
$\bar{P}$ is defined by turning the conditions \eqref{positivity} to weak inequalities, or taking the limits in $\mathbb{R}^m$ of sequences of $P$). Then
$\vol$ has well-defined extrema on $\bar{P}$, which are automatically critical if they belong to $P$. Because of the following proposition, the only
possibility for a critical point is to be an absolute maximum.

\begin{proposition}\label{prop:concavity}
 The volume of an ideal tetrahedron is a concave function of its dihedral angles. 
 \end{proposition}

\begin{proof} This follows from \fullref{tetvolume}, whose notations we use again: $x_t, y_t, z_t$ are the dihedral angles. By symmetry we may assume
$x_0,y_0\leq\pi/2$. Assume further that $x_t, y_t, z_t$ are affine functions of $t$ with first-degree coefficients $X,Y,Z$. \fullref{tetvolume}
implies $-d\vol/dt=X\log\sin x_t+Y\log\sin y_t+Z\log\sin z_t$, and by differentiating we obtain $$\eqalignbot{-d^2\vol/dt^2|_{t=0}&=X^2\cot
x_0+Y^2\cot y_0+Z^2\cot z_0\cr &=X^2\cot x_0+Y^2\cot y_0 +(X+Y)^2\frac{1-\cot x_0\cot y_0}{\cot x_0+\cot y_0}\cr &=\frac{(X+Y)^2+(X\cot
x_0-Y\cot y_0)^2}{\cot x_0+\cot y_0} \geq 0.\cr}\proved$$ \end{proof}

As a consequence, the volume functional $\vol$ is also concave on $\bar{P}$ and \fullref{main} holds if the maximum of $\vol$ is interior.
Next we explore the behavior of $\vol$ near the boundary of $\bar{P}$. 

\begin{proposition}[Simple degeneracy]\llabel{simpledeg} If $(Q_t)_{t\geq
0}$ is a smooth family of ideal tetrahedra with dihedral angles $x_t, y_t, z_t$ such that 
$x_0, y_0\in(0,\pi)$; $z_0=0$ and $\frac{dz_t}{dt}|_{t=0}>0$, then
$\frac{d\mathcal V}{dt}|_{t=0}=+\infty$. \end{proposition} 

\begin{proof} Simply check that the right hand side of \eqref{tetvolumeformula} goes to $0$ as $t$ goes to
$0$. We call this situation \emph{simple degeneracy\/} because the limiting triangle has only one vanishing angle (two of its vertices are therefore collapsed).
\end{proof} 

\begin{proposition}[Double degeneracy] \llabel{doubledeg} If $(Q_t)_{t\geq 0}$ is a smooth family of ideal tetrahedra satisfying $(x_0, y_0,
z_0)=(0,0,\pi)$ and $\frac{d}{dt}|_{t=0}(x_t, y_t, z_t)=(1+\lambda,1-\lambda,-2)$ with $\lambda \in (-1,1)$, then $$\exp \frac{-d\vol}{dt}\Big|_{t=0}=
\frac{1-\lambda^2}{4} \left ( \frac{1+\lambda}{1-\lambda} \right )^{\lambda}.$$ \end{proposition} 

\begin{proof} As $t$ goes to $0$, one has $\sin x_t \sim (1+\lambda)t$
and $\sin y_t \sim (1-\lambda)t$ and $\sin z_t\sim 2t$. The right hand side of \eqref{tetvolumeformula} is thus equivalent to $$\big((1+\lambda)t\big)^{1+\lambda}
\big((1-\lambda)t\big)^{1-\lambda} (2t)^{-2}=\frac{1-\lambda^2}{4} \left(\frac{1+\lambda}{1-\lambda} \right )^{\lambda}.$$ We call this situation \emph{double degeneracy\/}
because the limiting triangle has two vanishing angles (its vertices are distinct, but collinear). At a double degeneracy, the volume has directional derivatives, but no well-defined differential. \end{proof}

\section{Ruling out some degeneracies} \llabel{somemore}

From now on, we fix $(w_i)$ in the compact polyhedron $\bar{P}$ at a point which maximizes the total hyperbolic volume of all tetrahedra. To prove that $(w_i)$ is
critical for the volume $\vol$, we only need to make sure that $(w_i)$ lies in the interior $P$ of $\bar{P}$, ie that all $x_i, y_i, z_i$ lie in
$(0,\pi)$.

\begin{proposition} \llabel{tricot1} If for some index $i$, one of the numbers $x_i, y_i, z_i$ is $0$, then two of them are $0$ and the third is $\pi$. In other words,
there are no simple degeneracies, only double degeneracies. \end{proposition}
\begin{proof} If not, consider an affine segment from $(w_i)$ to some interior point of $P$. By \fullref{simpledeg}, the partial derivative of $\vol$
at $(w_i)$ along that segment is not bounded above, so $\vol$ was not maximal at $(w_i)$. \end{proof}

Tetrahedra $\Delta_i$ such that $(x_i,y_i,z_i)$ has one, and therefore two vanishing terms are called \emph{flat\/}, and are characterized by the property that
$w_i$ is either $0$ or $\frac{\pi}2$.
\begin{proposition}[Domino effect] \llabel{tricot2}  If two consecutive tetrahedra $\Delta_{i-1},\Delta_i$ are flat, then $\Delta_{i+1}$ is flat, too.
\end{proposition}
\begin{proof} We use only \fullref{interlettres} and the deductions recorded in \eqref{positivity}. There are three cases: \begin{enumerate} \item If $i$ is not a hinge index, flatness of $\Delta_i$
implies $w_{i-1}+w_{i+1}\in\{0,\pi\}$. By the range condition $0\leq w\leq \frac{\pi}2$, this implies $w_{i+1}\in\{0,\frac{\pi}2\}$, so $\Delta_{i+1}$ is flat. \item
If $i$ is a hinge index and $w_i=\frac{\pi}2$, we must have $|w_{i-1}-w_{i+1}|=\frac{\pi}2$, so by the range condition, $w_{i+1}$ is $0$ or $\frac{\pi}2$, and
$\Delta_{i+1}$ is flat. \item If $i$ is a hinge index and $w_i=0$, we have $|w_{i+1}-w_{i-1}|\leq 0$ so $w_{i+1}=w_{i-1}$. But $\Delta_{i-1}$ is assumed flat, and
therefore so is $\Delta_{i+1}$. Note that flatness of $\Delta_{i-1}$ was needed only in this case. \proved\end{enumerate}
\end{proof}
\begin{proposition} \llabel{tricot3} If $\Delta_i$ is flat, then $i$ is a hinge index and $w_i=0$. \end{proposition}
\begin{proof}In all other cases, the proof of \fullref{tricot2} actually forces $\Delta_{i+1}$ to be flat, which triggers a domino effect: all $\Delta_j$'s are
flat, and the volume is $0$ --- certainly not maximal. \end{proof}

\begin{trivial} \llabel{angsinedge} If  $ABC$ is a Euclidean triangle with positive angles and edge lengths $a,b,c$, then $$\begin{array}{ccccl} a=b &\iff
&\widehat{A}=\widehat{B} &\iff &\sin\widehat{A}=\sin\widehat{B} \\ a<b &\iff &\widehat{A}<\widehat{B} &\iff &\sin\widehat{A}<\sin\widehat{B}.\\ \end{array}$$
\end{trivial}

The volume $\vol$ is still supposed maximal, and we assume that some tetrahedra $\Delta_i$ are flat, ie that some hinge parameters $w_i$ vanish. Places where
a parameter $w_i$ vanishes will be signalled by a vertical bar: $\dots LL|RR\dots$ By \fullref{tricot2}, consecutive vertical bars are always separated by at
least two letters.

The patterns $RL|RL$ and $LR|LR$ can never occur, because increasing the incriminated $w_i$ to $\var$ would automatically increase the volume (note that
$w_{i-1}=w_{i+1}=:A$ by the hinge condition \eqref{positivity}-iii):
$$\begin{array}{c|ccccccccc}
\Omega& &L&             &R&    |    &L&             &R& \\
 w_i  &u& &      A      & & 0+\var  & &      A      & &v\\ \hline
 x_i  &.& &  u+A-\var   & &  \var   & &  \var+A-v   & &.\\
 y_i  &.& &\!\!-u+A+\var& &  \var   & &\!\!-\var+A+v& &.\\
 z_i  &.& &   \pi-2A    & &\pi-2\var& &   \pi-2A    & &.
\end{array}$$
This implies
$$\exp {\frac{-\partial \vol}{\partial \var}}\Big|_{\var=0} =\frac{1}{4}\cdot \frac{\sin(A-u)}{\sin(A+u)}\frac{\sin(A-v)}{\sin(A+v)}<1$$
where we used \fullref{angsinedge}, \fullref{tetvolume} and \fullref{doubledeg} (with $\lambda=0$).

Any vertical bar thus lives next to at least two consecutive identical letters (on at least one side). However, the patterns $R|LL|R$ and $L|RR|L$ are also prohibited
by \fullref{tricot3}, since the central (nonhinge) tetrahedron would have one vanishing angle ($a+c=0$ in the notations of \fullref{interlettres}).

\section{A geometrical lemma} \llabel{ageometricallemma}
\begin{definition} In the universal cover of the tessellation ${\mathcal A}$ of the torus at infinity of $V_{\varphi}$ (\fullref{pile}), a \emph{fan\/} is a sequence
of at least three consecutive layers, such that the first and last layers are grey and all layers in between are white. Fans are in bijection with the subwords of
$\Omega$ of the form $RL^kR$ or $LR^kL$ with $k\geq 2$ (see the remarks to \fullref{pile} in \fullref{combotorinf}). \end{definition}

\begin{lemma} \llabel{geolemma} Suppose $w_0=0$, so that $\Omega$ contains a subword $L|R^kL$ with $k\geq 2$, or $L|R^k|L$ with $k\geq 3$ (in the latter case, the
second bar indicates that $w_k$ vanishes as well as $w_0$). The corresponding fan admits a complete Euclidean structure with boundary (with angles prescribed by the
$w_i$'s). Moreover, let $Q,P,T$ be the lengths of the segments of the broken line corresponding to the first $R$, in the order indicated in \fullref{eventails} 
($P,T$ are the sides adjacent to the apex in a flat upward-pointing grey [hinge] triangle, in the sense of \fullref{fourtriangles}).
Then $Q<P+T$. \end{lemma}
\begin{figure}[ht!] \centering \includegraphics{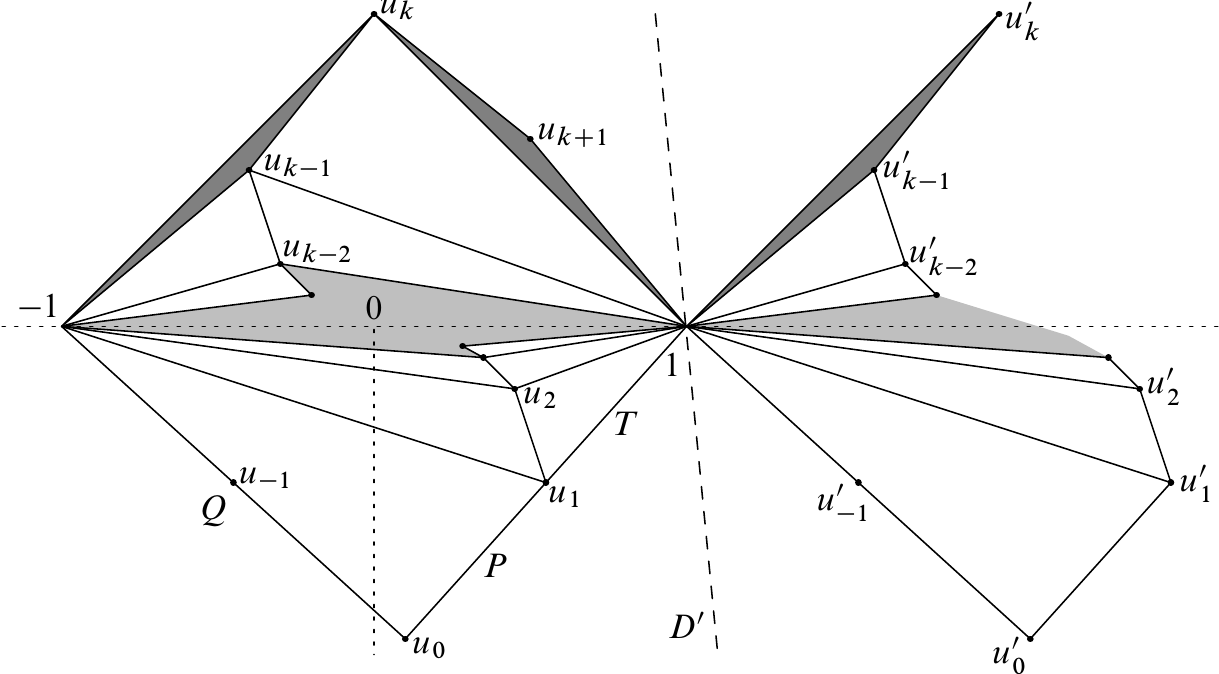} \caption{The situation where $Q\geq P+T$ cannot hold} \llabel{eventails}\end{figure}

\begin{proof} We first restrict our attention to the case $L|R^kL, ~k\geq 2$. The interior vertices of the topological fan correspond to the indices $i$ living between
two $R$'s, ie $1\leq i\leq k-1$, in the sense that holonomy around the $i$--th vertex $u_i$ (\fullref{eventails}) is controlled by $\partial \vol/\partial
w_i$ (\fullref{computhol}). When $2\leq i \leq k-1$, none of the triangles adjacent to $u_i$ are flat, so $w_i$ can vary in a small interval without making the
$m$--tuple $w$ exit the domain $\bar{P}$; consequently, the value of $w_i$ in that interval is critical, which by the first case (nonhinge) of \fullref{computhol}
implies that the holonomy around the associated vertex $u_i$ is trivial. As for $i=1$, the corresponding vertex $u_1$ is adjacent to a flat angle $z_0=\pi$ (\fullref{eventails}) so its holonomy is not imposed by the $w_i$'s (when a triangle has an angle $\pi$, the other two angles are always $0$ while the adjacent sides may
have arbitrary lengths). The case of $L|R^k|L$ is similar: for indices $2\leq i\leq k-2$, \fullref{computhol} applies, while for $i=1$ or $k-1$ holonomy is not imposed by the $w_i$'s.

Therefore we can embed the fan as an infinite necklace of triangulated polygons into ${\mathbb C}$. We shall no longer distinguish $L|R^kL$ from $L|R^k|L$ here, and shall formulate all properties in terms of complex numbers in order to keep track of \emph{oriented\/} angles. We make two consecutive nodes (ie lifts of the $(2k+4)$--valent
vertex of ${\mathcal A}$ associated to the full subword $LR^kL$) coincide with $-1$ and $1$ in $\mathbb{C}$, and also denote by $u_i$ the complex coordinate of the copy of $u_i$ between these nodes (the index $i$ actually ranges from $-1$ to $k+1$; see \fullref{eventails}). Incidentally, one can show that the $u_i$ form (part of) an orbit of a certain deck transformation of the universal cover $\mathbb{H}^3\rightarrow V_{\varphi}$.  We will discuss this more at the end of \fullref{numericalexample}.

We arrange matters so that ${\rm Im}(u_0)<0$ and $u_1$ lies on the open
segment $(u_0,1)$. While removing any node disconnects the fan, Condition \eqref{coherence}-iii implies that the image of the embedding is invariant under a horizontal
translation of length 2. In particular, the geometric link of each node, such as $1$ or $-1$, is completely determined. To prove the assertion of the Lemma, it is
sufficient to show that ${\rm Re}(u_0)<0$. Assume ${\rm Re}(u_0)\geq 0$ (so $u_0$ lies in the lower-right quadrant) and aim for a contradiction.

The similarity property of the triangles with vertices $1,u_i,u_{i+1}$ and $-1,u_i,u_{i-1}$ is expressed by the relation
$\frac{u_{i+1}-1}{u_i-1}=\frac{u_{i-1}+1}{u_i+1}$, hence by induction $$(u_{i+1}-1)(u_i +1)=(u_i-1)(u_{i-1}+1)= \dots=(u_1-1)(u_0+1)=:K.$$ Then, as $u_1$ sits between
$1$ and $u_0$, the number $u_1-1$ is a positive (real) multiple of $u_0-1$, so $K$ is a positive multiple of $u_0^2-1$ which implies ${\rm Im}(K)\leq 0$. Let $D$ be
the line through $0$ and the points $\pm \sqrt{K}$: either $D$ is vertical, or $D$ visits the upper-left and lower-right (open) quadrants. Let $D'$ be the line through
$1$, parallel to $D$; and define $u'_i:=2+u_i$ for all $i$. By definition of $K$, the rays $[1,u_{i+1})$ and $[1,u'_i)$ are symmetric with respect to $D'$. Moreover,
the rays through $u_0,u_1,\dots,u_k,u_{k+1},u'_k,u'_{k-1},\dots,u'_0,u'_{-1},u_0$ issuing from $1$ (in that cyclic order) divide ${\mathbb C}$ into angular sectors of
sum $2\pi$ realizing the geometric link of a node of the fan, as specified by the angles $x_i, y_i, z_i\geq 0$. Finally, since all these angles are nonnegative, the
symmetry of the link with respect to $D'$ implies that for all $-1\leq i\leq k$, the point $u'_i$ (resp.\ $u_{i+1}$) is on the right (resp.\ left) of $D'$.

Recall ${\rm Im}(u_0)<0$. Suppose by induction ${\rm Im}(u_i)<0$ for some $0\leq i \leq k$. Then ${\rm Im}(u'_i)<0$. Considering the direction of $D'$ and the symmetry property with
respect to $D'$, this implies ${\rm Im}(u_{i+1})<0$. By an immediate induction, the angular sector $\smash{\widehat{u_{k+1} 1 u'_k}}$ (just above $1$) is larger than
$\pi$. But it is an angle of the link at the node $1$ (namely, $z_{k+1}$), giving a contradiction.\end{proof}

Of course, a statement similar to \fullref{geolemma} holds for subwords $LR^k|L$, and also for $R|L^kR, R|L^k|R$ and $RL^k|R$.

\section{Ruling out all degeneracies} \llabel{hallali}
\begin{trivial} \llabel{lambda} If $U$ and $V$ are positive constants, the function defined on $(-1,1)$ by
$$f(\lambda):=\frac{1-\lambda^2}{4}U(\frac{1+\lambda}{1-\lambda}V) ^{\lambda}$$ takes the value $\frac{U}{(1+V)(1+V^{-1})}$ for some $\lambda$. It is in fact an absolute minimum:
indeed, $$(\log f)'(\lambda)=\log \Big[\smash{\frac{1+\lambda}{1-\lambda} V}\Big],$$ so $f$ is minimal when the bracket is $1$, and the result follows by direct computation.
\end{trivial}

Now we can prove that the configuration $\dots RR|L\dots$ (and similarly $\dots LL|R\dots$) never occurs, which will imply \fullref{main}. The strategy is to
suppose $RR|L$ appears, ie $w_j=0$ for some $j$ (for notational convenience we assume $j=2$). Next, replace $w_2$ by $\var$ and $w_1$ by $w_1+\lambda\var$, for a wisely
chosen $\lambda$. The volume $\vol$ will then increase.
(The value of $\lambda$, which does not need to be explicitly computed, will maximize $\partial \vol / \partial \var$, and $e^{-\partial\vol/\partial\var}$ will be the
value of $f$ given in \fullref{lambda}.  We will specify in due time what the parameters $U,V$ are.) 
Volume computations follow from \fullref{doubledeg} (at the index $i=2$) and
\fullref{tetvolume} (other indices).

\subsection[Case 1: RR|LR]{Case 1: $RR|LR$} \llabel{case1} According to \fullref{interlettres}, the angles are as follows. Note the relation $w_1=w_3=:A$, a consequence of the hinge condition
\eqref{positivity}-iii.
$$\begin{array}{c|ccccccccc}
\Omega&                  &R&                       &R&       |       &L&         &R& \\
     i&         0        & &           1           & &       2       & &    3    & &4\\
   w_i&         u        & &     A+\lambda\var     & &    0+\var     & &    A    & &v\\ \hline
   x_i& \xi-A-\lambda\var& &-u+2A+2\lambda\var-\var& &(1-\lambda)\var& & \var+A-v& &.\\
   y_i&\eta+A+\lambda\var& &         u+\var        & &(1+\lambda)\var& &-\var+A+v& &.\\
   z_i&      \pi-2u      & &  \pi-2A-2\lambda\var  & &  \pi-2\var    & &  \pi-2A & &.\end{array}$$
 We have $u=w_0>0$ because $L|RR|LR$ is impossible (according to the discussion after \fullref{angsinedge}). Thus, $\Delta_0$ is not flat: if $-1<\lambda<1$, then $\varepsilon$ can take small positive values without making any of the $x_i, y_i, z_i$ negative. With the correct choice of $\lambda$, we deduce
\begin{eqnarray*} \exp \frac{-\partial \vol}{\partial \var}\Big|_{\var=0} &=& \frac{1-\lambda^2}{4}
\overbrace{\underbrace{\frac{\sin(A-v)}{\sin(A+v)}}_{\leq 1} \underbrace{\frac{\sin y_1}{\sin x_1}}_{Q/P}}^U \left ( \frac{1+\lambda}{1-\lambda}
\underbrace{\frac{\sin y_0 \sin^2 x_1}{\sin x_0 \sin^2 z_1}}_{P/T=:V} \right )^{\lambda} \\ &\leq & \frac{Q/P}{(1+P/T)(1+T/P)} = \frac{1}{1+P/T}~\frac {Q}{P+T} < 1
\end{eqnarray*} by \fullref{geolemma} (see \fullref{finale}\,(left) --- again, the sine relation in triangles was used to compute $Q/P$ and $P/T$). 
\begin{figure}[ht!] \centering \includegraphics{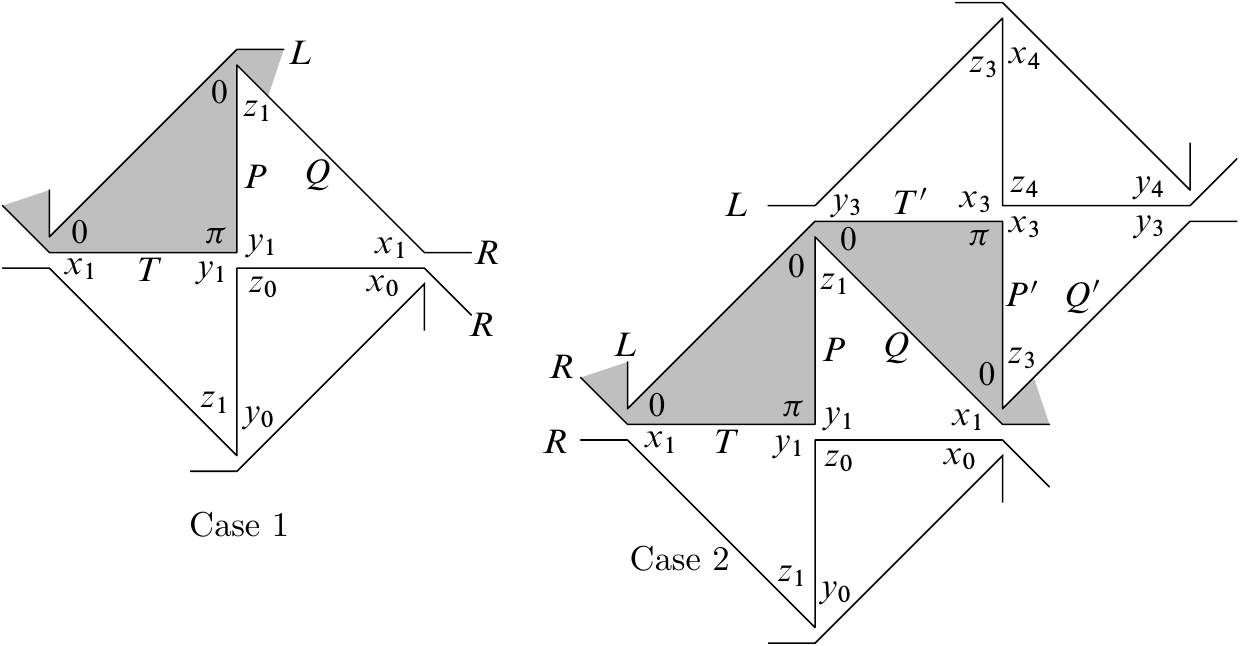} \caption{}\llabel{finale} \end{figure}

\subsection[Case 2: RR|LL]{Case 2: $RR|LL$} \llabel{case2}

$$\begin{array}{c|ccccccccc}
\Omega&             &\!\!\!\!R\!&             &\!\!\!\!R\!&  |  &\!\!\!\!L\!& &\!\!\!\!L\!&  \\
    i &         0        & &           1           & &       2       & &     3    & &   4    \\
  w_i &         u        & &    A+\lambda\var      & &    0+\var     & &     A    & &   v    \\ \hline
  x_i &\xi-A-\lambda\var & &-u+2A+2\lambda\var-\var& &(1-\lambda)\var& &  \var+v  & & A+\xi' \\
  y_i &\eta+A+\lambda\var& &        u+\var         & &(1+\lambda)\var& &-\var+2A-v& &-A+\eta'\\
  z_i &      \pi-2u      & &  \pi-2A-2\lambda\var  & &  \pi-2\var    & &  \pi-2A  & & \pi-2v \end{array}$$
 First consider the value of $A$. Since there can be no vertical bars immediately before or after $RR|LL$, the tetrahedra $\Delta_0,\Delta_1,\Delta_3,\Delta_4$
have positive angles, so the parameter $A$ (which does not contribute to the angles of any other tetrahedra) can vary freely in an open interval when $\var=0$. So $A$
must be critical giving $$ 1 = \exp \frac{-\partial \vol}{\partial A}\Big|_{\var=0} =\underbrace{\frac{\sin y_0\sin^2 x_1}{\sin x_0\sin^2 z_1}}_{P/T}
\underbrace{\frac{\sin^2 y_3\sin x_4}{\sin^2 z_3\sin y_4}}_{P'/T'}.$$ Hence $P/T=T'/P'$. Therefore, with the right choice of $\lambda$, \begin{eqnarray*} \exp
\frac{-\partial \vol}{\partial \var}\Big|_{\var=0} &=& \frac{1-\lambda^2}{4} \overbrace{\underbrace{\frac{\sin y_1}{\sin x_1}}_{Q/P} \underbrace{\frac{\sin x_3}{\sin
y_3}}_{Q'/P'}}^{U} \left ( \frac{1+\lambda}{1-\lambda} \underbrace{\frac{\sin y_0 \sin^2 x_1}{\sin x_0 \sin^2 z_1}}_{P/T=T'/P'=:V} \right )^{\lambda} \\ & = &
\frac{Q/P~~\cdot ~~Q'/P'}{(1+T/P)(1+T'/P')} = \frac{Q}{P+T} ~\frac{Q'}{P'+T'} < 1 \end{eqnarray*} by \fullref{geolemma}; see \fullref{finale}\,(right). 

We conclude with two remarks. First, up to replacing the monodromy $\varphi$ with $\varphi^2$ (thus doubling the period $m$ of $\Omega$), we can always assume $\Omega$ has at least $6$ letters: that way, all columns of the tables above are neatly distinct, and to recover the original $V_{\varphi}$ we just quotient out by the extra symmetry (which the volume maximizer $(w_i)_{i\in\mathbb{Z}/m\mathbb{Z}}$ must respect, by concavity of the volume functional $\mathcal{V}$). Compare with the remark that closes \fullref{positiveangles}. Second, the choice of $\lambda \in (-1,1)$, which may seem ``magical'' at first glance, is essentially our only degree of freedom in \fullref{case1}--\fullref{case2}: the volume is already assumed critical with respect to most parameters (including $A$, the common value of $w_1$ and $w_3$); therefore, only deformations of $w_1,w_2,w_3$ need to be considered, and if we assume $\smash{\frac{\partial w_2}{\partial\varepsilon}}=1$, then only the value of the difference $\frac{\partial w_1}{\partial \varepsilon} - \frac{\partial w_3}{\partial \varepsilon}\in (-1,1)$ matters.

\fullref{main} is proved.

\section[A numerical example]{A numerical example: $R^NL^M$} \llabel{numericalexample}

In this section we fix two large enough integers $N$ and $M$ and investigate the behavior of the angles for $\Omega=R^NL^M$: the angles made positive by the previous
computations will turn out to be very small. We will directly construct a Euclidean realization of \fullref{pile}, automatically unique up to isometry. Since $N$
and $M$ are large, there exist small complex numbers $a,a',b,b'$ such that
\begin{equation} \llabel{plugin} \left \{ \begin{array}{lcr} \sin a &=& i \tan b \cos b' \\ \sin a' &=& -i \tan b' \cos b \end{array} \right .
\quad \text{where} \quad \left \{ \begin{array}{lcl} b &=& (\pi-2a)/N \\ b' &=& (\pi-2a')/M. \end{array} \right . \end{equation}
A way to compute $a,a'$ is to set $a_0=a'_0=0$ and to define inductively $$a_{s+1}=\arcsin\Big(i\tan\frac{\pi-2a_s}{N}\cos\frac{\pi-2a'_s}{M}\Big)$$ and a similar expression for
$a'_{s+1}$, with a change of sign. The sequences $a_s, a'_s$ converge exponentially fast to $a,a'$. The constants $a,a'$ become arbitrarily small for large enough $N,M$, hence $b\sim\pi/N$
and $b'\sim \pi/M$. So plugging into \eqref{plugin}, $a\sim i\pi/N$ and $a'\sim -i\pi/M$. Using the Landau symbol $O(A,B)$ in the sense of $O(\max\{A,B\})$, this in turn
yields
\begin{equation} \llabel{estime} \begin{array}{rcl}
b\, =&\frac{  \pi}{N}-\frac{2i\pi}{N^2}&+O(\frac{1}{N^3},\frac{1}{M^3}) \\ \\
b'  =&\frac{  \pi}{M}+\frac{2i\pi}{M^2}&+O(\frac{1}{N^3},\frac{1}{M^3}) \end{array} \qquad \begin{array}{rcl}
a\, =&\frac{ i\pi}{N}+\frac{2 \pi}{N^2}&+O(\frac{1}{N^3},\frac{1}{M^3}) \\ \\
a'  =&\frac{-i\pi}{M}+\frac{2 \pi}{M^2}&+O(\frac{1}{N^3},\frac{1}{M^3}).
\end{array} \end{equation}
(In fact $a,a',b,b'$ are analytic functions of $\frac{1}{N},\frac{1}{M}$, by the Implicit Function Theorem.)
\begin{proposition} The fan which corresponds to $R^N$ can be embedded into ${\mathbb C}$ with nodes at complex coordinates $\pm \cot b$ and intermediary vertices
$\cot(a+sb)$ where $-1\leq s \leq N+1$; similarly, the fan corresponding to $L^M$ can be embedded into ${\mathbb C}$ (possibly with a different scaling factor) with nodes $\pm \cot b'$ and intermediary vertices
$\cot (a'+sb')$ where $-1\leq s \leq M+1$ (see \fullref{twofans}). \end{proposition}
\begin{figure}[ht!] \centering \includegraphics{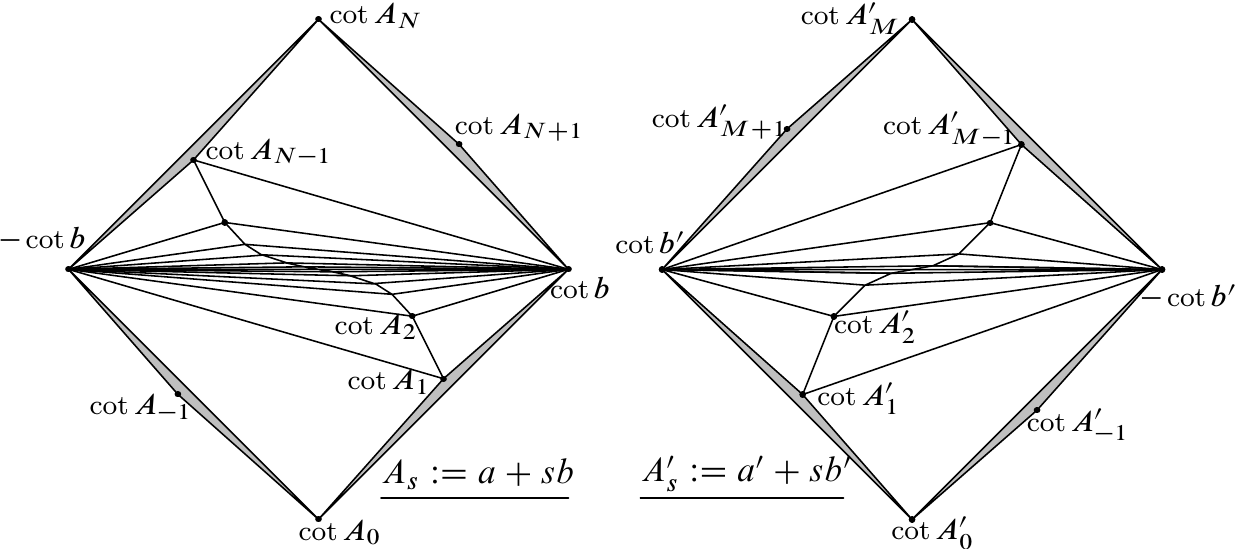} \caption{}\llabel{twofans}\end{figure}

\begin{proof} There are several things to check. First, the congruence of pairs of triangles inside each fan follows from the identity of complex ratios
\begin{equation} \llabel{squaredsines} \frac{\cot(A_s+b)-\cot b}{\cot(A_s+b) -\cot A_s}=\frac{\cot(A_s-b)+\cot b}{\cot(A_s-b)-\cot A_s}=\frac{\sin^2 A_s}{\sin^2 b} \end{equation}
where $A_s=a+sb$ for $0\leq s \leq N$, and an identical relation for $a', b'$.

Next, each fan, when stripped of two of its four limiting (grey) triangles, is a parallelogram. This follows from $\cot(a+Nb)=\cot(\pi-a)=-\cot a$, and the same for
$a'$. In particular, each fan admits a center of symmetry.

Furthermore, these two parallelograms are congruent. To see this, let $\alpha, \alpha', \beta, \beta'$ denote the \emph{squared cotangents\/} of $a,a',b,b'$. Raising
\eqref{plugin} to the power $-2$, we get $$ \left \{ \begin{array}{lcl} 1+\alpha&=&-\beta(1+\beta'^{-1}) \\ 1+\alpha'&=&-\beta'(1+\beta^{-1}) \end{array} \right .
\quad\text{hence}\quad \frac{\alpha}{\alpha'}=\frac{\beta+\beta/\beta'+1}{\beta'+\beta'/\beta+1} = \frac{\beta}{\beta'}, $$ so $\unfrac{\cot a}{\cot
b}=\pm\unfrac{\cot a'}{\cot b'}$, the correct sign being minus by the estimates \eqref{estime}.

Further yet, the limiting (grey) triangles of the two fans have the same shape: by \eqref{squaredsines} their complex ratios are $\unfrac{\sin^2 a}{\sin^2 b}$ and
$\unfrac{\sin^2 b'}{\sin^2 a'}$, both equal by \eqref{plugin} to $-\unfrac{\cos^2 b'}{\cos^2 b}$.

Finally, all triangles are positively oriented, that is,  ${\rm Im}(\unfrac{\sin^2 A_s}{\sin^2 b})>0$ for all $0\leq s\leq N$. We first check this for $s=0$ (the case $s=N$ will follow by symmetry): we have $\unfrac{\sin^2 a}{\sin^2 b}=-\punfrac{\cos^2 b'}{\cos^2 b}=-\punfrac{(1-\sin^2 b')}{(1-\sin^2 b)}$. In the latter expression, both
the numerator and denominator are $\sim 1$, but their imaginary parts are equivalent to $-4\pi^2/M^3$ and $4\pi^2/N^3$ respectively (with an $O(N^{-4},M^{-4})$ error), by
\eqref{estime}. Therefore, the ratio does lie above the real line, and the ``pinched'' angles of the limiting (grey) triangles in \fullref{twofans} are both
roughly
$$2\pi^2(N^{-3}+M^{-3}) \text{ radians}.$$  Very pinched, but not flat!
\begin{figure}[ht!]
\centering \includegraphics{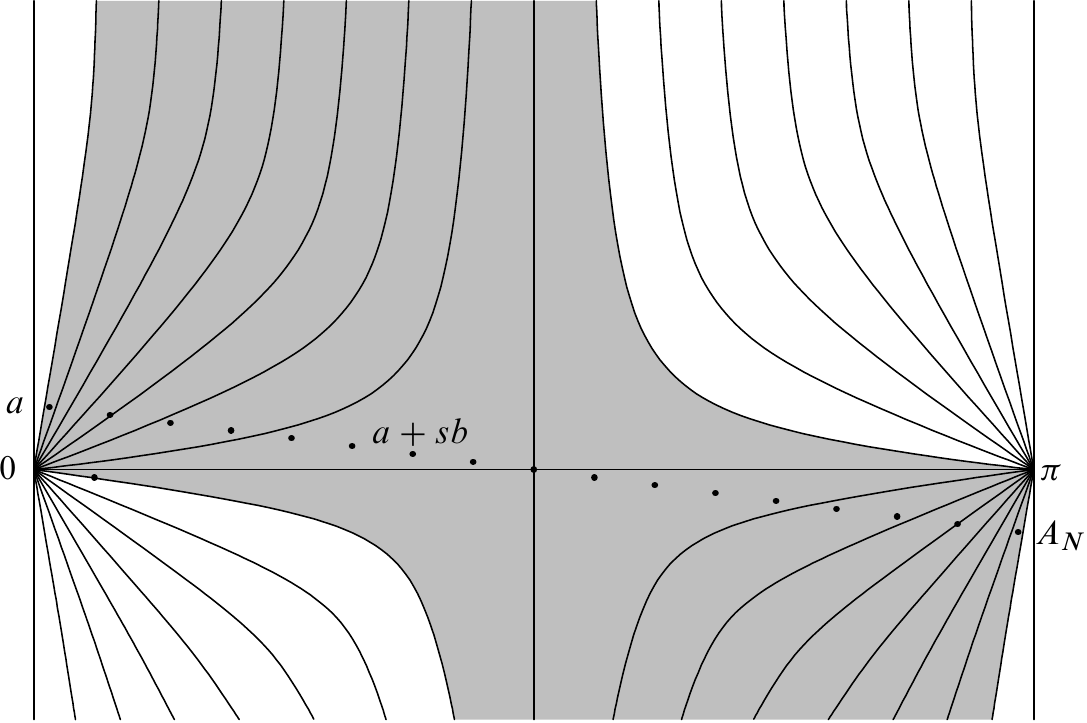}  \caption{Level curves of
$z\mapsto \arg(\sin z)$} \llabel{levelcurves} \end{figure}

To check that $\arg(\unfrac{\sin A_s}{\sin b})$ lies in $(0,\pi/2)$ for all other $0\leq s \leq N$, we need to draw the level curves of $z\mapsto \arg (\sin z)$ in
${\mathbb C}$. This is done in \fullref{levelcurves}, in the case $0\leq {\rm Re}\,(z) \leq \pi$: the curves fall into $4$ symmetric families (``quadrants'' meeting at $\pi/2
\in {\mathbb C}$), and it is an easy exercise to check that the families above (resp.\ below) the real axis are made of convex (resp.\ concave) curves. The authorized
region for the $A_s=a+sb$ is in grey (with a narrow collar near $\pi/2$); the forbidden regions are left in white. The segment $[A_0A_N]$ clearly stays in the grey
region, which implies the result.

By the same argument, the triangles in the fan of $L^M$ are well-oriented, too. Therefore we may tile the plane with parallelograms (or fans) congruent to those in
\fullref{twofans} to get a Euclidean realization of \fullref{pile}. \end{proof}

Finally, notice that the Kleinian group associated to the embedding of the left fan in \fullref{twofans} contains the M\"obius transformation $z\mapsto\punfrac{z\cos
b-\sin b}{(z\sin b+\cos b)}$ (it sends each tetrahedron sitting above a triangle in the left half of the fan to the tetrahedron sitting above the similar triangle in the
right half). Therefore, $2ib$ (and similarly $-2ib'$) are the complex lengths of very short closed geodesics in the hyperbolic manifold $V_{\varphi}$.

A consequence of \fullref{main} is that the volume of any angle structure (defined by \fullref{interlettres}, where \eqref{positivity} holds) is a lower bound for the volume of the manifold $V_{\varphi}$. One gets bounds which are sharp in terms of distances in the Farey graph; see \fullref{volumebundles}.

\section{Once-punctured tori and $4$--punctured spheres} \llabel{spheres}

\fullref{main} is still true if we replace the once-punctured torus $T$ by the $4$--punctured sphere $S$, and the map $\varphi\co T \rightarrow T$ by an
orientation-preserving homeomorphism $\varphi_S\co S \rightarrow S$ (of course, we must specify how the ``eigenvalues'' of $\varphi_S$ are defined). In fact, the tetrahedra
of the resulting manifold $V_{\varphi_S}$ and of $V_{\varphi}$ are metrically the same; only the combinatorics of their gluing changes a little.

Define $R:={\mathbb R}^2 \setminus {\mathbb Z}^2$ and the maps $\alpha,\beta,\sigma \co R \rightarrow R$ characterized by $\alpha(x)=x+(1,0)$; $\beta(x)=x+(0,1)$;
$\sigma(x)=-x$. Then we have natural identifications $T=R/\langle\alpha,\beta \rangle$ and $S=R/\langle\alpha^2,\beta^2,\sigma \rangle$. Define also $T':=R/\langle
\alpha^2,\beta^2 \rangle$ (note that $T'$ is a $4$--punctured torus).

One can show that up to isotopy, any orientation-preserving diffeomorphism $\varphi_S$ of $S$ lifts to a map $\varphi_R\co R \rightarrow R$ such that
$\varphi_R(x)=Mx+v$ for some $M\in \SL_2({\mathbb Z})$ and $v\in {\mathbb Z}^2$. Moreover, $(\pm M)$ and $(v \text{ mod } 2{\mathbb Z}^2)$ are unique. So we may define
the eigenvalues of $\varphi_S$ (up to sign) as those of $M$. Observe finally that $\varphi_R$ induces orientation-preserving diffeomorphisms $\varphi'\co T' \rightarrow
T'$ and $\varphi_T\co T \rightarrow T$. 

There are obvious coverings $T'\rightarrow T$ and $T'\rightarrow S$, of degrees $4$ and $2$ respectively. Given an ideal triangulation $\tau_T$ of $T$ (corresponding
to a Farey triangle), we can lift $\tau_T$ to an ideal triangulation $\tau'$ of $T'$. Observe that $\sigma$ acts on $T'$ as a properly discontinuous involution fixing
$\tau'$, hence $\tau'$ descends to a triangulation $\tau_S$ of $S$. It is easy to see that $\tau_S$ has the combinatorics of a tetrahedron. If $\tau_T$ and $\tau^1_T$
are separated by a diagonal exchange (see \fullref{leucocephale} for a definition), then $\tau_S$ and the corresponding $\tau^1_S$ are separated by two diagonal
exchanges on opposite edges. Mutatis mutandis, the construction of \fullref{topotorinf} provides ideal (topological) triangulations of $T' \times {\mathbb
R}$ and $S \times {\mathbb R}$, as well as of their quotients $V_{\varphi'}$ and $V_{\varphi_S}$. There are coverings $V_{\varphi'}\rightarrow V_{\varphi_S}$ and
$V_{\varphi'}\rightarrow V_{\varphi_T}$ (of degrees $2$ and $4$ respectively), and all these manifolds are hyperbolic when $\varphi_S$ has distinct real eigenvalues.


\appendix\def\leftmark{David Futer}
\def\rightmark{Canonical triangulations (appendix)}
\section{Geometric triangulations of two-bridge link\\ complements}
  \llabel{appendix}

\cl{\sc David Futer}

This appendix applies Gu\'{e}ritaud's techniques to find geometric  
triangulations for the hyperbolic two-bridge knot and link  
complements. These ideal triangulations are, in essence, the  
monodromy triangulations of $4$--punctured sphere
bundles, closed off in a slightly different way. They were  
constructed and studied in great detail by Sakuma and Weeks \cite 
{sakuma-weeks}. Akiyoshi, Sakuma, Wada and Yamashita have announced  
a proof that these triangulations are, in fact, geometrically  
canonical \cite{aswy-announce}, which is a stronger statement than  
our result.

We will begin by reviewing two-bridge links and these ideal  
triangulations. We will then explain how the methods of the preceding  
paper give linear angle structures for these
triangulations and prove that the volume function is maximized in the  
interior of the space of angle structures. Finally, we will prove two  
corollaries of this argument: a two-sided bound on the volume of the  
link complement and a result about arcs in the plane being hyperbolic  
geodesics.

\subsection{Braids and two-bridge links} Let $S$ be a $4$--punctured  
sphere, visualized concretely as a square pillowcase with the corners  
removed. A $4$--string braid between two
pillowcases, one interior and one exterior, defines a so-called \emph 
{product region} $S \cross I$. We will restrict our attention to  
alternating braids in which the top
right strand is free of crossings. (See \fullref{fig:pillow-braid}\,(a).)

\begin{figure}[ht!]
\labellist
\small\hair 2pt
\pinlabel (a) [br] at 99 468
\pinlabel (b) [br] at 236 468
\endlabellist
\begin{center}
\includegraphics{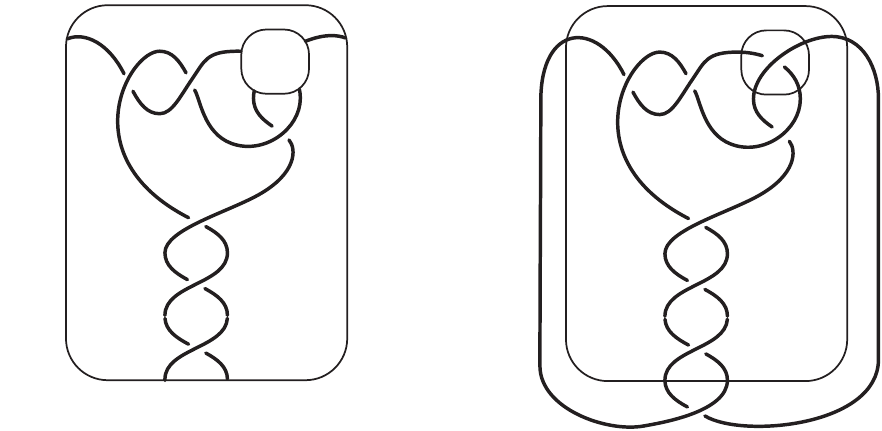}
\end{center}
\caption{(a) An alternating braid between two pillowcases,  described  
by the word $\Omega = R^3 L^2 R$ (b) The corresponding two-bridge  
link $K(\Omega)$}
\llabel{fig:pillow-braid}
\end{figure}

The mapping class $\varphi\co S \to S$ induced by this braid can be  
described by a word
$$
\Omega := \left\{ \begin{array}{llll}
R^{a_1} L^{a_2} \cdots R^{a_n} & \textrm{ or } & L^{a_1} R^{a_2}  
\cdots L^{a_n} & \textrm{ with odd $n$, or} \\
R^{a_1} L^{a_2} \cdots L^{a_n} & \textrm{ or } & L^{a_1} R^{a_2}  
\cdots R^{a_n} & \textrm{ with even $n$.} \\
\end{array}\right.
$$
Here, $R$ encodes a crossing on the bottom pair of strands, and $L$  
encodes a crossing on the left pair of strands, as in \fullref 
{fig:rl-crossings}. Each \emph{syllable\/}
of $\Omega$ (that is, each maximal subword $R^{a_i}$ or $L^{a_i}$) 
corresponds to a \emph{twist region\/} in which two strands of the  
braid wrap around each other $a_i$ times.
As we read $\Omega$ from left to right, we scan the crossings from  
the outside in. Note that, unlike the case of punctured torus  
bundles, our word $\Omega$ has a beginning
and an end. For concreteness, we will focus on the case when $\Omega$  
starts with $R$, as in \fullref{fig:pillow-braid}\,(a); the $L$--case is similar.

\begin{figure}[ht!]
\labellist
\small\hair 2pt
\pinlabel $R$ at 277 656
\pinlabel $L$ at 406 656
\endlabellist
\begin{center}
\includegraphics{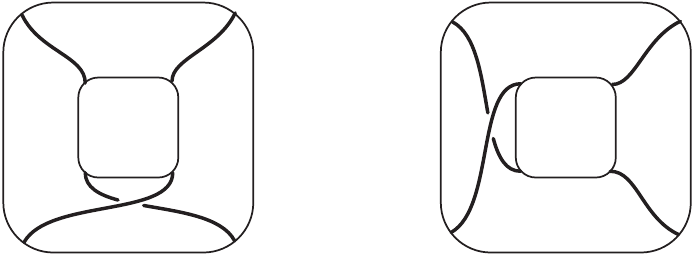}
\caption{The letters $R$ and $L$ acting on strands of a braid}
\llabel{fig:rl-crossings}
\end{center}
\end{figure}

An alternating braid of this sort can be completed to a link diagram,  
as follows. Outside the outer pillowcase, we connect the bottom left  
strand to the top right, and the
bottom right strand to the top left, adding an extra crossing. (Up to  
isotopy of $\mathbb{S}^2$, there is a unique way to do this while keeping an  
alternating projection. In \fullref{fig:pillow-braid}\,(b), the extra crossing was arbitrarily placed  
at the bottom of the diagram.) Similarly, inside the inner pillowcase  
we connect the strands in a
diagonal fashion, adding an extra crossing while preserving the  
alternating projection. This creates an alternating link $K(\Omega)$,  
as in \fullref{fig:pillow-braid}\,(b).

$K$ is called a \emph{two-bridge link\/} because this diagram can be  
isotoped so that the pillowcases are horizontal, and the connecting  
strands form two bridges between the
strands of the braid. It is well-known that, apart from the trivial  
link of one or two components, every two-bridge link can be  
constructed in this way. (See, for example,
Murasugi \cite[Theorems 9.3.1 and 9.3.2]{murasugi:knot-book}.)

William Menasco's theorem about hyperbolic alternating links \cite 
{menasco-alt} contains the following special case.

\begin{theorem}\llabel{thm:hyperbolic-links}
The two-bridge link $K(\Omega)$ is hyperbolic if and only if $\Omega 
$ has two or more syllables.
\end{theorem}

Just as with punctured torus bundles, we will give a direct proof of  
the ``if'' direction of this theorem by finding a geometric ideal  
triangulation of the link complement.
The ``only if'' direction is immediate: a word with a single syllable  
produces a link with a single twist region, which must be a torus link.

\subsection{The ideal triangulation}
The word $\Omega$ describes a monodromy triangulation of the product  
region $S \cross I$, in exactly the same fashion as for $4$--punctured  
sphere bundles. In fact, because this
product region is the complement of a braid in the part of $\mathbb{S}^3$ bounded by the two pillowcases,
we can locate the edges of the triangulation concretely in the projection diagram.

Let $c = \sum_{i=1}^n a_i$ be the length of $\Omega$. Each letter $ 
\Omega_i$ ($1 \leq i \leq c$), and thus each crossing in the  
alternating braid, corresponds to a
$4$--punctured sphere $S_i \subset \mathbb{S}^3 \setminus K$, with the four  
strands of $K$ seen in \fullref{fig:rl-crossings} passing through  
the four punctures. The braid
induces an ideal triangulation on each $S_i$, whose edges come from  
arcs in the diagram that look vertical and horizontal immediately before
and/or after the corresponding crossing. See \fullref{fig:near-crossing}.

\begin{figure}[ht!]
\begin{center}
\includegraphics{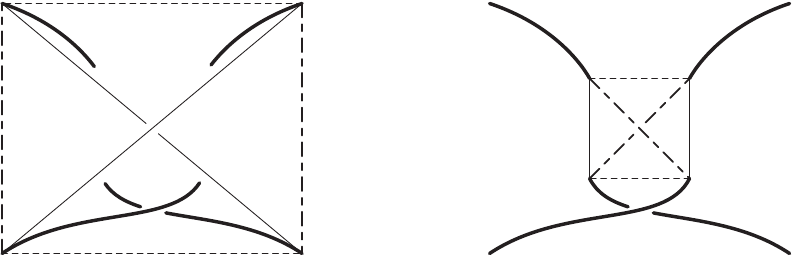}
\end{center}
\caption{Two views of the same $4$--punctured sphere $S_i$ living near  
a crossing in the link diagram (arcs with the same  
dashing pattern have the same slope)}
\llabel{fig:near-crossing}
\end{figure}

Just as with punctured torus bundles, we can locate the progressively  
changing triangulations in the Farey tesselation $F$. Each  
triangulation of a $4$--punctured sphere $S_i$,
containing six edges of three different slopes, corresponds to a  
Farey triangle $t_i$. Triangles $t_i$ and $t_{i+1}$ share an edge $e_i 
$, whose endpoints are the shared
slopes of $S_i$ and $S_{i+1}$. If $\Omega_i = R$, the path from $e_ 
{i-1}$ to $e_i$ takes a right turn across $t_i$; if $\Omega_i = L$,  
the path takes a left turn. This rule
also defines an \emph{initial edge\/} $e_0$, because $\Omega_1 = R$, so  
a right turn should take $e_0$ to $e_1$. Similarly, the action of $ 
\Omega_c$ defines a \emph{terminal
edge} $e_c$.

For each $e_i$, $1 \leq i \leq c-1$, the $4$--punctured spheres $S_i$  
and $S_{i+1}$ are joined together along four edges, two for each  
endpoint of $e_i$. In between them lies a
layer $\Delta_i = \Delta(e_i)$ of two tetrahedra, whose bottom  
surface $S_i$ has the triangulation of $t_i$ and whose top surface $S_ 
{i+1}$ has the triangulation of
$t_{i+1}$. (See \fullref{fig:tetrahedron-layer} for an example.)  
Stacking these tetrahedron layers together produces a layered  
triangulation of the product region between
$S_1$ and $S_c$.

\begin{figure}[ht!]
\begin{center}
\includegraphics{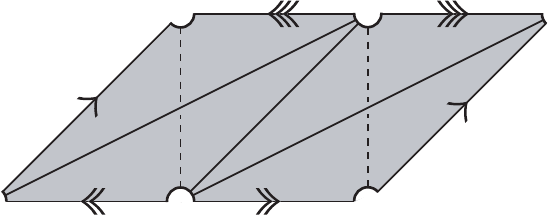}
\caption{The tetrahedron layer $\Delta_1 = \Delta(e_1)$, made of two  
tetrahedra contained between $4$--punctured spheres $S_1$ and $S_2$  
(sides with identical arrows are
identified)}
\llabel{fig:tetrahedron-layer}
\end{center}
\end{figure}

If we wanted to construct a bundle of $4$--punctured spheres over the  
circle, we would glue the top of this product region to the bottom.  
To recover the complement of the
link $K(\Omega)$, we follow a slightly different procedure. On the $4$--punctured sphere $S_1$, corresponding to the first crossing inside  
the outer pillowcase, let the
\emph{peripheral edges\/} be the edges whose slope is the vertex of $t_1 
$ opposite the initial edge $e_0$.

We will fold $S_1$ along the two peripheral edges, identifying its  
ideal triangles in pairs. \fullref{fig:fold-clasp} shows that this  
creates exactly the desired effect of
connecting the strands of $K$ in pairs, with a twist. This full twist  
corresponds to the first two crossings in the link projection: the  
first crossing in the braid, as well
as the ``extra'' crossing outside the outer pillowcase. The four nonperipheral edges on $S_1$ are identified to a single edge, isotopic  
to a short arc near the crossing. The
two ideal triangles that remain after folding are clasped together  
around this edge, which we call the \emph{core\/} of the clasp.

\begin{figure}[ht!]
\labellist
\small\hair 2pt
\pinlabel $\Rightarrow$ at 135 197
\pinlabel $\Downarrow \quad \text{isotopy}$  at 218 126
\pinlabel $\Leftarrow$ at 135 54
\endlabellist
\begin{center}
\includegraphics{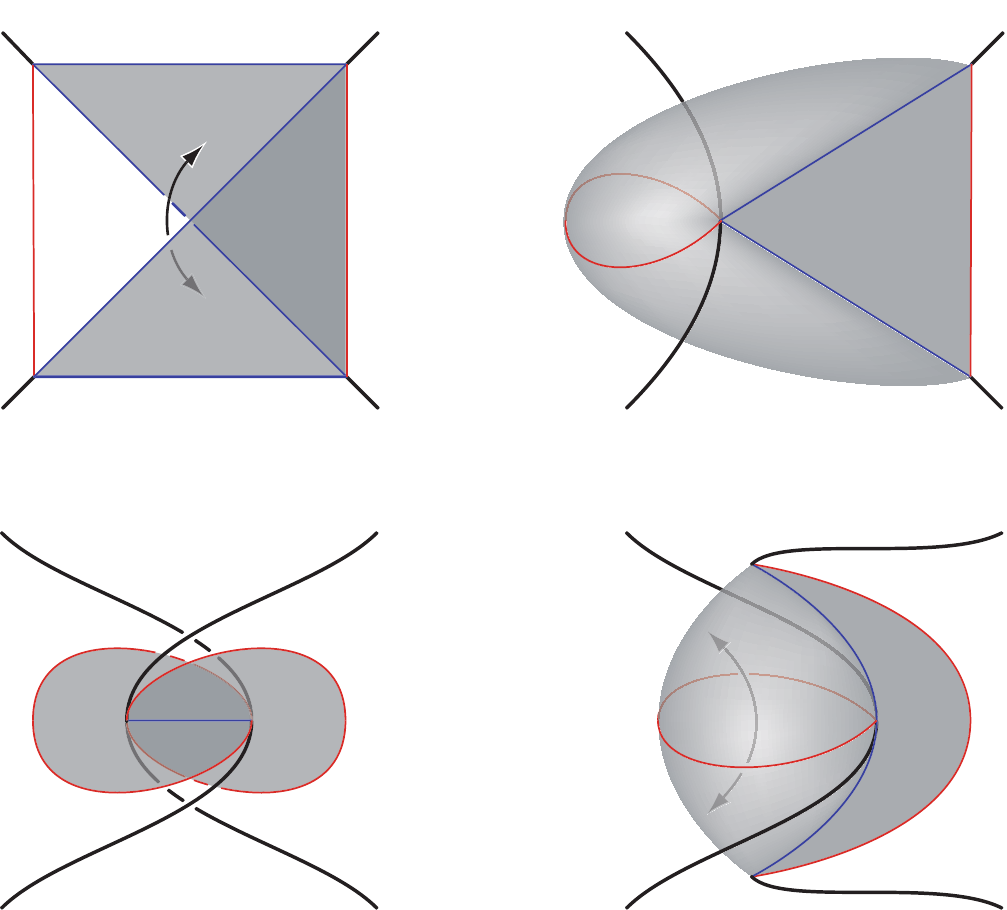}
\caption{Folding the pleated surface $S_1$ produces the first two  
crossings in the link.}
\llabel{fig:fold-clasp}
\end{center}
\end{figure}

Let $\alpha$ be the mirror image of the peripheral slope across $e_0$  
in the Farey graph. Topologically, folding $S_1$ as above amounts to  
attaching a thickened disk of
boundary slope $\alpha$ to the outer pillowcase \cite[Lemma II.2.5] 
{sakuma-weeks}.

In a similar fashion, we define the peripheral edges on $S_c$ to be  
the edges whose slope is  the vertex of $t_c$ opposite $e_c$. We fold  
$S_c$ along these two peripheral
edges, identifying its faces to a clasp of two ideal triangles.  
Topologically, attaching $2$--handles to $S_1$ and $S_c$ results in a  
space homeomorphic to the link complement.
Combinatorially, folding $S_1$ and $S_c$ defines a gluing pattern for  
all the faces of the tetrahedra, giving us an ideal triangulation of  
$\mathbb{S}^3 \setminus K$. See
\cite[Section II.2]{sakuma-weeks} for more details of this  
triangulation.

\subsection{Combinatorics at the cusp}

To describe the combinatorics of the boundary component(s) of $\mathbb{S}^3  
\setminus K$, we will first focus on the product region between the  
pleated surfaces $S_1$ and $S_c$. In
the layered triangulation of this product region, each layer $\Delta_i 
$ consists of two tetrahedra, $D_i$ and $D_i'$, as in \fullref 
{fig:tetrahedron-layer}. It is clear
from the figure that each tetrahedron has exactly one vertex at each  
puncture of $S_i$, ie at each strand of the $4$--string braid between  
$S_1$ and $S_c$. Since the
combinatorics of the four strands are identical, let us focus on a  
single puncture of the $4$--punctured sphere.

The tetrahedron layer $\Delta_i$ intersects the neighborhood of a  
puncture in two \emph{boundary triangles\/}, one from a truncated  
vertex of $D_i$ and one from $D_i'$. These
boundary triangles meet at two vertices that come from shared edges  
of $D_i$ and $D_i'$. (This completes a loop, corresponding to the  
meridian of a component of $K$.) The
apices of the two triangles point in different directions, as in  
\fullref{fourtriangles}. As with  
punctured torus bundles, these layers of
boundary triangles are stacked together, forming fans that correspond  
to syllables of the word $\Omega$.

\begin{figure}[ht!]
\begin{center}
\includegraphics{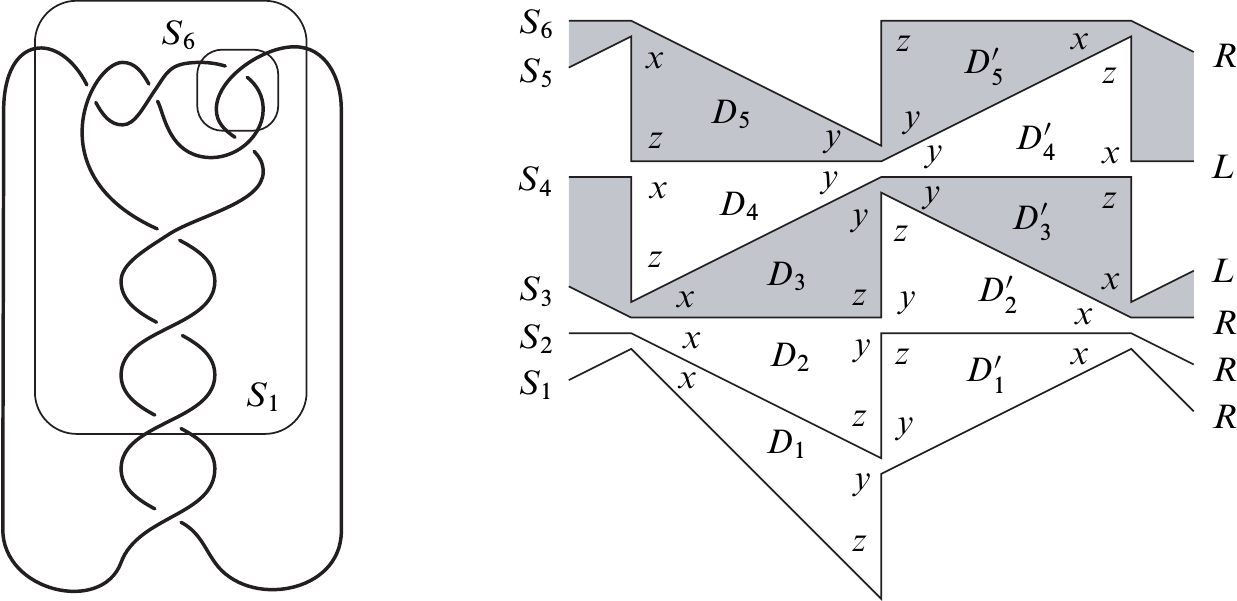}
\caption{Left:\ \ the link $K(\Omega)$ corresponding to $\Omega = R^3  
L^2 R$\newline Right:\ \ a cusp view of the product region between $S_1$ and  
$S_c$}
\llabel{fig:cusp-braid}
\end{center}
\end{figure}

The resulting cusp triangulation which corresponds to the product region  
of the link $K(R^3 L^2 R)$ is shown in \fullref{fig:cusp-braid}.  
As in \fullref{pile}, the
triangles are shown opened up, and the hinge layers are shaded.  
Observe that this picture is combinatorially equivalent to the 
corresponding picture for punctured torus bundles, quotiented by the
hyperelliptic involution.

We have labeled the dihedral angles of the tetrahedra of $\Delta_i$  
by numbers (``angles'') $x_i$, $y_i$, $z_i$, following the same  
conventions as in \fullref{pile}. Note that
our choices of dihedral angles force the tetrahedra $D_i$ and $D_i'$  
to be isometric; this does not impede the goal of finding a geometric  
triangulation.  In the sequel, we
will not distinguish between $D_i$ and $D_i'$.

\begin{figure}[ht!]
\begin{center}
\includegraphics[width=4.5in]{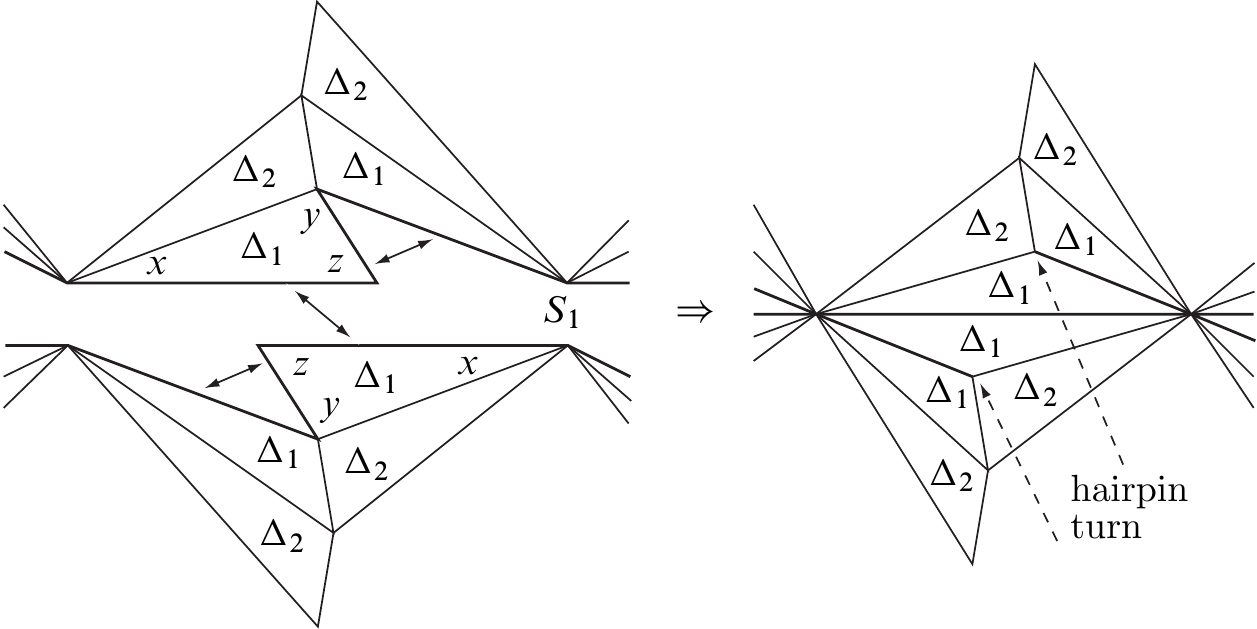}
\caption{A cusp view of the folding that occurs at a clasp}
\llabel{fig:cusp-clasp}
\end{center}
\end{figure}

When the pleated surface $S_1$ is folded to form a clasp, the zigzag  
line in which it intersects the cusp also becomes folded, creating a  
``hairpin turn.'' Because this
folding procedure joins the punctures of $S_1$ in pairs, as in \fullref{fig:fold-clasp}, the boundary triangles on those punctures are  
also joined together. The resulting
cusp triangulation in the neighborhood of $S_1$ can be seen in \fullref{fig:cusp-clasp}. At the other end of the product region, the  
clasp of $S_c$ appears on the cusp in
the same way.

Once we have folded the clasp surfaces as prescribed, the truncated  
vertices of the tetrahedra combine to form either a single torus that  
traverses the product region four
times (in case $K$ is a knot), or two tori, each of which traverses  
the product region twice (in case $K$ is a two-component link). In  
either case, the local combinatorics are
the same, and the affine equations that the dihedral angles of the  
tetrahedra must satisfy are derived in the same way.

To find the hyperbolic structure on $\mathbb{S}^3 \setminus K$, we will study  
the space $P$ of \emph{angle structures\/} for the triangulation, ie  
the space of positive dihedral
angles that line up correctly around each edge. We need to complete  
three steps:
\begin{enumerate}
\item Parameterize $P$ and check that it is nonempty, as in \fullref{positiveangles}.

\item Show that a critical point of the volume functional on $P$  
gives a complete hyperbolic metric, as in \fullref 
{hyperbolicvolume}.

\item Prove that at any point on the boundary of $\bar{P}$, the  
volume can be increased by unflattening the flat tetrahedra, as in  
\fullref{somemore}--\fullref{hallali}.
\end{enumerate}

This will imply that the volume functional is maximized in the  
interior of $P$, guaranteeing a critical point and thus a hyperbolic  
metric. For each of the three steps, the
argument is essentially the same as Gu\'{e}ritaud's.

\subsection{Angle structures and volume}

Following \fullref{positiveangles} of the main article, we will  
parameterize the dihedral angles of the tetrahedra by \emph{pleating  
angles} on the pleated $4$--punctured
spheres. Each sphere $S_i$ described above has a natural transverse  
orientation that points (equivalently) toward the inside of the link  
projection, toward increasing
indices, and upward in \fullref{fig:cusp-braid}. Just as in  
\fullref{positiveangles}, for any edge $e \subset S_i$, we can  
define the \emph{pleating angle\/} $\alpha$ to
be the (signed) external angle at $e$, with signs chosen so that $ 
\alpha$ is positive whenever the angle above $e$ is less than $\pi$.

On the pleated sphere $S_{i+1}$ living between $\Delta_i$ and $\Delta_ 
{i+1}$, this definition will give pleating angles  $$-2w_i, \quad  
2w_{i+1} \quad \mbox{and}
\quad 2w_i - 2w_{i+1},$$ exactly as in \eqref{pleatingangles}.  
The clasp surface $S_1$, which borders $\Delta_1$ on one side and is  
folded on the other side, will have
pleating angles $-\pi$, $2w_1$, and $\pi - 2w_1$, where the pleating  
angle of $-\pi$ corresponds to the hairpin turn in \fullref 
{fig:cusp-clasp}. Thus, if we define $w_0 =
\pi/2$ (even though there is no tetrahedron layer $\Delta_0$), the  
pleating angles on $S_1$ will be given by the same expressions as  
above. Similarly, setting $w_c = \pi/2$
allows us to label the pleating angles on $S_c$ by the same  
expressions as in \eqref{pleatingangles}.

\begin{lemma}\llabel{lemma:angle-struct-param}
For $i=0, \ldots, c$, choose a parameter $w_i$, such that $w_0 = w_c  
= \pi/2$ and $w_1, \ldots, w_{c-1}$ satisfy the range, concavity, and  
hinge conditions \eqref{positivity}.
For each such choice of parameters, set the dihedral angles of the  
tetrahedra as in \fullref{interlettres}. Then
\begin{enumerate}
\item for each $\Delta_i$, the angles $x_i$, $y_i$, $z_i$ are  
positive and add up to $\pi$,
\item the dihedral angles around each edge add up to $2\pi$, and
\item for each $S_i$, the pleating angles add up to $0$.
\end{enumerate}
\end{lemma}

\begin{proof}
The range, concavity, and hinge conditions imply that all the  
tetrahedron angles $x_i$, $y_i$, $z_i$ are positive, and the claim  
that their sum is $\pi$ is immediate from
\fullref{interlettres}. For each pleated surface $S_i$, the  
pleating angles sum to $0$ by construction. Therefore, it remains to  
check the angle sum around each edge $e$ of
$\mathbb{S}^3 \setminus K$.

If $e$ is not the core of $S_1$ or $S_c$, the combinatorial picture  
of \fullref{fig:cusp-braid} is the same as the one for torus  
bundles. Thus, as in \fullref{positiveangles}, each layer $\Delta_i$ contributes precisely the  
difference between the pleating angles of the neighboring surfaces,  
and the sum around $e$ simplifies to
$2\pi$. If $e$ is the core of a clasp, say the core of $S_1$, the  
left panel of \fullref{fig:cusp-clasp} shows that four sectors  
contribute dihedral angles at $e$: two
sectors that have angle $z_1$, plus two fans of angles above pleated  
surface $S_1$, each having dihedral angle $2w_1$. Thus, because $z_1  
+ 2w_1 = \pi$, the total angle sum
at $e$ is $2\pi$. \end{proof}

By \fullref{lemma:angle-struct-param}, our triangulation will have  
an angle structure whenever we set $w_0 = w_c = \pi/2$ and  
interpolate between these parameters in a way
that satisfies the range, concavity, and hinge conditions. One can  
always do this graphically, by first fixing $w_i$ for the hinge  
indices and then connecting the hinges by
pieces of parabolas, as in \fullref{dubya} (this is where we use  
the hypothesis that $\Omega$ contains at least one hinge).

\begin{lemma}\llabel{lemma:twobridge-critical}
Let $P$ be the open affine polyhedron of angle structures for the  
triangulation of $\mathbb{S}^3 \setminus K$, parameterized by sequences
$(\frac{\pi}{2},w_1,\dots,w_{c-1},\frac{\pi}{2})$, as in \fullref 
{lemma:angle-struct-param}. Then a point of $P$ is a critical point  
of the volume functional $\vol$ if and only if the corresponding
tetrahedron shapes give a complete hyperbolic structure on $\mathbb{S}^3  
\setminus K$.
\end{lemma}

Just like \fullref{rivinlemma}, this is an instance of a much more  
general theorem of Rivin, Chan, and Hodgson \cite{chan-hodgson,rivin}. It is also possible to prove
\fullref{lemma:twobridge-critical} directly, using the same line of  
argument as in \fullref{rivinlemma}, although in the setting of  
two-bridge links this would require
considering a number of special cases.

\subsection{Flat tetrahedra never maximize volume}

The proof that the maximum of $\vol$ occurs in the interior of $P$  
closely tracks \fullref{somemore}--\fullref{hallali} of Gu\'{e}ritaud's paper. We begin by ruling out many
types of degeneracies on $ \bdy \bar{P}$.

\begin{lemma}\llabel{lemma:classify-degeneracy}
Let $(w_1, \ldots, w_{c-1})$ be the point of  $\bar{P}$ at which  
the volume functional $\vol$ attains its maximum. Then $(w_i)$ has  
the following properties:
\begin{enumerate}
\item For each $i$, if one of $x_i$, $y_i$, $z_i$ is $0$, then two  
are $0$, ie $\Delta_i$ is \emph{flat\/}.
\item If $\Delta_i$ is flat, then $w_i = 0$.
\item If $\Delta_i$ is flat, then $i$ is a hinge index, not equal to  
$1$ or $c-1$.
\item If $\Delta_i$ is flat, then $i$ is adjacent to at least two  
consecutive identical letters.
\item If both hinge layers at the ends of a syllable $R^k$ or $L^k$  
are flat, then $k \geq 3$.
\end{enumerate}
\end{lemma}

\begin{proof}
All the discussion and results of \fullref{somemore} apply  
equally well in our context. Thus we have conclusion $(1)$ as a  
restatement of \fullref{tricot1}. The
domino effect of \fullref{tricot2} also applies; in fact, it  
is clear from the proof of the Proposition that the domino effect  
works both forward and backward. Thus it
does not matter that our word $\Omega$ is not cyclic.

Almost all the claims of $(2)$--$(5)$ are proved in \fullref 
{somemore}, either in \fullref{tricot3} or in the discussion  
that follows. The one exception is the claim
that $\Delta_1$ and $\Delta_{c-1}$ cannot flatten. This follows  
because, in the case of two-bridge links, $w_0 = w_c = \pi/2$. Thus  
setting $w_1$ or $w_{c-1}$ to $0$ or
$\pi/2$ will trigger the domino effect and force all the tetrahedra  
to flatten, giving a volume of $0$. \end{proof}

It remains to prove that at any point $(w_i) \in \bdy \bar{P}$  
satisfying the properties of \fullref{lemma:classify-degeneracy},  
the volume will increase as we move
into the interior of $P$. Following Gu\'{e}ritaud, we do this using  
the geometrical statement of \fullref{geolemma}. The proof of this  
lemma transfers perfectly to our
context when the fan under consideration corresponds to a subword in  
the interior of $\Omega$. As it turns out, the same statement is even  
easier to prove when the degenerate
layer is near the beginning or end of $\Omega$.

\begin{lemma}\llabel{lemma:last-case}
Recall the word $\Omega = R^{a_1} L^{a_2} \cdots$, and suppose that  
the hinge layer $\Delta_{a_1}$ has flattened, with $w_{a_1} = 0$.  
Then the fan corresponding to $R^{a_1}L$
admits a complete Euclidean structure with boundary along $S_{a_1}$.  
Let $Q$, $P$, $T$ be the lengths of the segments of the broken line  
in which $S_{a_1}$ intersects the
cusp, as in \fullref{fig:terminal-fan}. Then $Q < P+ T$.
\end{lemma}

Of course, the analogous statement holds near the end of $\Omega$.

\begin{figure}[ht!]
\begin{center}
\includegraphics[width=3.7in]{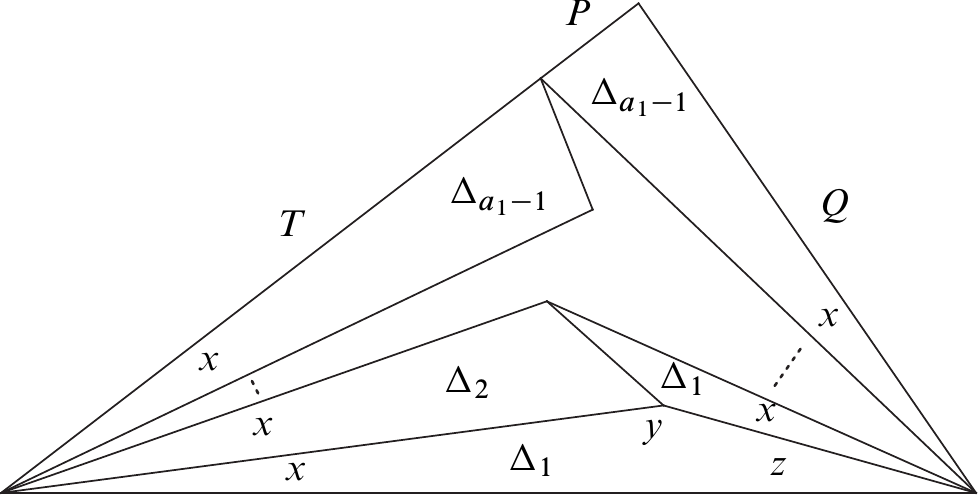}
\caption{In a fan at the beginning of $\Omega$, $Q < P+T$.}
\llabel{fig:terminal-fan}
\end{center}
\end{figure}

\begin{proof}
Note that, by \fullref{lemma:classify-degeneracy}, we have $a_1  
\geq 2$. For all $1 \leq i \leq a_1-2$, the parameter $w_i$ can vary  
in a small interval, and no
generality is lost in assuming that $\partial \vol / \partial w_i = 0 
$. (Otherwise, volume is easy to increase.)  Then, as in the proof of  
\fullref{geolemma}, the
criticality of volume with respect to these parameters implies that  
the fan of $R^{a_1}L$ has a complete Euclidean structure. If $i>1$,  
the computation is the same as
in \fullref{computhol}. If $i=1$, the angle $z_0$ of \fullref 
{holonomies} is replaced by the hairpin turn, and the factor $\frac 
{\sin y_0 }{\sin x_0}$
disappears from the computation of $\exp(-\partial\vol/\partial w_1) 
$. Thus we once again have $1=\exp(-\partial\vol/\partial w_1)=|\psi 
(\alpha)|$, where $\psi$ is the
reduced holonomy and $\alpha$ is a loop around the vertex at the  
hairpin turn.

In our situation, the fan of $R^{a_1}L$ is a (tessellated) Euclidean  
triangle, as in \fullref{fig:terminal-fan}, in which $Q$ and $P+T$ are two of the  
sidelengths. In the triangular fan, the angle opposite $Q$ is $x_1 +  
\ldots + x_{a_1 - 1}$, while the angle
opposite $P+T$ is $z_1 + x_1 + \ldots + x_{a_1 - 1}$. Thus, because  
its side is opposite the smaller angle, $Q < P+ T$. \end{proof}

Armed with \fullref{geolemma} and its analogue in \fullref 
{lemma:last-case}, we can complete the proof of \fullref 
{thm:hyperbolic-links} by following the argument of
\fullref{hallali}. For each hinge index $j$ with $w_j = 0$, we  
pick a vector along which to deform the neighboring parameters, in a  
way that will maximize the derivative
of volume. Deforming in this optimal direction with speed $\varepsilon 
$, we can compute $\partial \vol/\partial \varepsilon$ in terms of  
the geometry of the fan(s) that
adjoin $\Delta_j$.

In most situations, the exponentiated derivative appears as a product  
of the exact same sine ratios as in \fullref{hallali}. The only  
exception occurs when $j=2$
or $j=c-2$, because there are no tetrahedra corresponding to $w_0$ or  
$w_c$. If (without loss of generality) $j=2$, the angle $z_0$ is  
replaced by the hairpin turn, while
the factor $\frac{\sin y_0}{\sin x_0}$ disappears from the  
computation of $P/T$, just as in the proof of \fullref{lemma:last-case}. Thus $\partial \vol/\partial
\varepsilon$ has the same expression in terms of $Q,Q',P,P',T,T'$ as  
in \fullref{hallali}.

In every case, the geometric statement $Q < P+ T$, applied to two  
fans if necessary, implies that $\partial \vol/\partial \varepsilon >  
0$. Therefore $\vol$ is maximized at a
critical point where all tetrahedron angles are positive, completing  
the proof of \fullref{thm:hyperbolic-links}.

\section{Applications}
Our construction of geometric triangulations using volume  
maximization methods has several corollaries that relate hyperbolic  
geometry (of bundles and links) to combinatorics (of the Farey  
complex and link diagrams). 

\subsection{Volume estimates for bundles}
\llabel{volumebundles}

\begin{theorem}\label{cor:bundle-volume}
Let $V_\varphi$ be a once-punctured torus bundle defined by the cyclic  
word $\Omega = R^{a_1}L^{b_1} \cdots R^{a_n}L^{b_n}$. Then
$$2 n\, v_3 \; \leq \;  \mathrm{Vol}(V_\varphi) \;<\; 2 n \, v_8,$$
where $v_3 \approx 1.0149$ is the volume of a regular ideal  
tetrahedron and $v_8 \approx 3.6638$ is the volume of a regular ideal  
octahedron. Both of these bounds are sharp.
\end{theorem}

This is a sharp, quantitative version of Brock's result in \cite 
{brock}, in the special case of punctured torus bundles. The upper bound is not new; 
it is a special case of \cite[Corollary 2.4]{agol}.

\begin{proof}
To prove the lower bound on volume, we exhibit a particular angle  
structure. It follows from the proof of \fullref{main} that the  
volume of the complete hyperbolic structure on $V_\varphi$ is the  
global maximum of $\vol$ over the closed polytope $\bar{P}$. Thus, for  
any $w \in \bar{P}$, $\vol(w) \leq  \mathrm{Vol}(V_\varphi)$.

We choose $w$ in the simplest possible manner: by letting $w_i = \frac 
{\pi}{3}$ for all $i$. It is easy to check that this choice of  
parameters satisfies at least the weak form of all the inequalities  
of \eqref{positivity}. In other words, all angles are nonnegative,  
and $w \in \bar{P}$. Any nonhinge tetrahedron $\Delta_i$ will have one  
vanishing angle and will therefore have volume $0$. However, every  
hinge tetrahedron $\Delta_i$ will be a regular tetrahedron, with all  
angles $\frac{\pi}{3}$ and volume $v_3$. Thus $\vol(w) = 2nv_3$.

When $\Omega = (RL)^n$, and thus all tetrahedra are hinges,  
this choice of angles will in fact give the complete structure on $V_ 
\varphi$. These bundles, which are $n$--fold cyclic covers of the  
figure--8 knot complement, have volume exactly $2nv_3$. On the other  
hand, the bundles that contain some nonhinge tetrahedra have $\mathrm 
{Vol}(V_\varphi) > 2nv_3$, because $\vol$ is maximized at a point  
where all tetrahedra have positive angles.

To prove the upper bound, we employ a Dehn surgery construction. For  
every Farey vertex $s$ which corresponds to a fan $RL^*R$ or $LR^*L$   
(with $*>0$), drill a closed curve of slope $s$ out of the fiber at  
the corresponding level. (These are exactly the closed curves whose  
length was estimated in \fullref{numericalexample}.) The  
resulting \emph{drilled bundle\/} turns out to be an $n$--fold cyclic  
cover of the Borromean rings complement, with volume $2 n v_8$.

We can recover $V_\varphi$ by Dehn filling the extra cusps of the  
drilled bundle. Because volume goes down under Dehn filling, $\mathrm 
{Vol}(V_\varphi) < 2 n v_8$. If we pick a word $\Omega$ with very  
long syllables, as in \fullref{numericalexample}, the resulting  
bundle $V_\varphi$ has volume arbitrarily close to the volume of this  
surgery parent.
\end{proof}

\begin{corollary}\label{cor:sphere-bundle-volume}
Let $V_\varphi$ be a $4$--punctured sphere bundle defined by the cyclic  
word $\Omega = R^{a_1}L^{b_1} \cdots R^{a_n}L^{b_n}$. Then
$$4 n\, v_3 \; \leq \;  \mathrm{Vol}(V_\varphi) \;<\; 4 n \, v_8.$$
\end{corollary}

\subsection{Volume estimates for links}

The volumes of link complements can also be estimated in terms of  
diagrams. We say that a link diagram $D$ is \emph{reduced\/} if no  
single crossing separates $D$. Its \emph{twist number\/} $\mathrm{tw}(D) 
$ is the number of equivalence classes of crossings, where two  
crossings are considered equivalent if there is a loop in the  
projection plane intersecting $D$ transversely precisely in the two  
crossings. When the diagram $D$ depicts a two-bridge link  
constructed from a braid, as in \fullref{fig:pillow-braid}, $\tw(D) 
$ is precisely the number of syllables of the word $\Omega$  
describing the braid.

\begin{theorem}\label{thm:2bridge-volume}
Let $D$ be a reduced alternating diagram of a hyperbolic two-bridge  
link $K$. Then
$$2v_3 \, \tw(D) - 2.7066 \; < \;  \mathrm{Vol}(\mathbb{S}^3 \setminus K) \;< 
\; 2 v_8 (\tw(D)-1).$$
The upper bound is sharp, and the lower bound is asymptotically sharp.
\end{theorem}

There are known diagrammatic volume bounds for the general  
class of alternating links. On the lower side, Agol,  
Storm, and W\,Thurston proved that the volume of an alternating link  
is at least $\frac{v_8}{2}(\tw(D)-2)$ \cite{ast-guts}, improving an earlier bound due to Lackenby \cite{lack-volume}. On the upper side, Agol and D\,Thurston  
proved an asymptotically sharp bound of $10v_3   
(\tw(D)-1)$ \cite[Appendix]{lack-volume}. Numerically, \fullref{thm:2bridge-volume} improves the multiplicative constant in the lower bound from 1.8312 to 2.0299 and the multiplicative constant in the upper bound from 10.1494 to 7.3277 (in the special case of $2$--bridge links).

Compared to its predecessors, \fullref{thm:2bridge-volume} is less general, but uses only very elementary methods. The proof in \cite{ast-guts} relies in a fundamental way on Perelman's results about the  
monotonicity of volume under Ricci flow with surgery. By contrast,  
the lower bound in \fullref{thm:2bridge-volume} only relies on  
the explicit study of angled triangulations.

\begin{figure}[ht!]
\begin{center}
\includegraphics{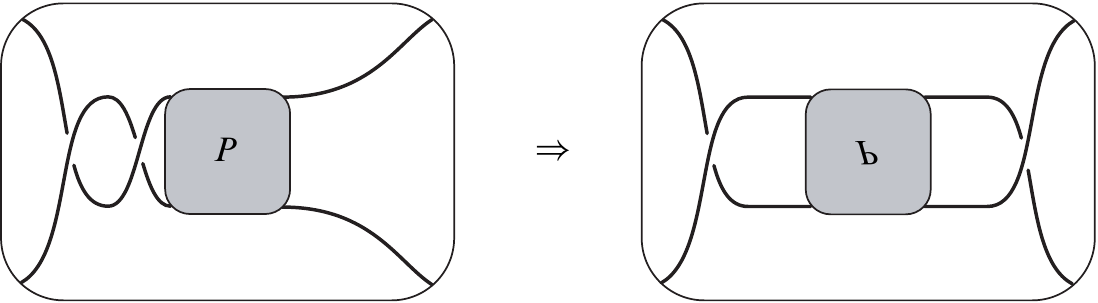}
\caption{A flype in a diagram of a two-bridge link. The shaded  
rectangle $P$ is a pillowcase of the braid.}
\llabel{fig:flype}
\end{center}
\end{figure}

\begin{proof}
First, we claim that it suffices to consider a diagram $D$ like the  
one in \fullref{fig:pillow-braid}, constructed from an alternating  
braid with one free strand.  By the Menasco--Thistlethwaite flyping  
theorem \cite{mt:flyping}, any pair of reduced alternating diagrams  
of $K$ are related by a sequence of \emph{flypes\/}, as in \fullref 
{fig:flype}. It is easy to check that the twist number of a diagram  
is invariant under flypes. Thus the number of syllables of the word $ 
\Omega$ is the twist number of \emph{any\/} reduced alternating diagram  
of $K$.

\medskip {\bf Lower bound}\qua
Suppose that the link is defined by the word $\Omega = R^{a_1} L^ 
{a_2} \cdots R^{a_n}$ (the parity of $n=\tw(D)$ and the letter of the  
first and last syllables are unimportant). As in the proof of \fullref{cor:bundle-volume}, we explicitly choose a point $w = (\frac{\pi} 
{2}, w_1, \ldots, w_{c-1}, \frac{\pi}{2})$ of $\smash{\bar{P}}$. However, the  
concavity condition of \eqref{positivity} does not allow us to set  
$w_i=\frac{\pi}{3}$ when $i$ is too close to $0$ or $c$. Instead, we  
proceed as follows. We let $w_i=\frac{\pi}{3}$ for all $a_1 \leq i  
\leq c-a_n$. For the indices of the first and last fans, we  
interpolate linearly between $\frac{\pi}{2}$ and $\frac{\pi}{3}$. As  
before, it is easy to check that these values of the parameters make  
all tetrahedron angles nonnegative, and give us a point of $\smash{\bar{P}}$.

When $i$ is a hinge and $a_1 < i < c-a_n$, the two tetrahedra of $ 
\Delta_i$ have all dihedral angles $\frac{\pi}{3}$, and volume $v_3$.  
For $n \geq 3$, there are exactly $n-3$ hinge indices of this type.  
When $i=a_1$ or $i=c-a_n$, we can compute from \fullref 
{interlettres} that the three angles of $\Delta_i$ are
$$\frac{\pi}{3} + t_i, \quad \frac{\pi}{3} - t_i \quad \mbox{and}  
\quad \frac{\pi}{3}, \quad \mbox{where} \quad
t_i=|w_{i+1} - w_{i-1}| \leq \frac{\pi}{6}.$$
By \fullref{prop:concavity}, the volume defined by these  
angles is smallest at the extreme value of $t_i=\frac{\pi}{6}$, when  
the three angles are $\frac{\pi}{2}, \frac{\pi}{3}, \frac{\pi}{6}$.  
Still assuming $n \geq 3$, by Formula \eqref{untetraedre} the four  
tetrahedra in the two terminal hinge layers each have volume at least  
$0.84578$... Putting it all together gives
$$\vol(w) \; > \; 2v_3 \, (\tw(D)-3) + 4 \times 0.84578 \; > \; 2v_3  
\, \tw(D) - 2.7066.$$
(As a special case, if $n=2$, $\vol(w) > 2 \times 0.84578$ also  
satisfies the theorem.)

To prove that this bound is asymptotically sharp, let $\Omega=(RL)^k 
$, for large $k$. Then $\tw(D)=2k$, and the triangulation of $K 
(\Omega)$ consists of $2(\tw(D)-1)$ tetrahedra, all of them hinges.  
Since the volume of an ideal tetrahedron is bounded above by $v_3$, $ 
\mathrm{Vol}(\mathbb{S}^3 \setminus K) \leq 2v_3 (\tw(D)-1)$, a value
whose ratio to the lower bound of the theorem approaches $1$ as $\tw 
(D)$ gets large.

\begin{figure}[ht!]
\labellist
\small\hair 2pt
\pinlabel $K$ [tl] at 171 692
\pinlabel $J$ [tl] at 282 692
\pinlabel $L$ [tl] at 391 692
\endlabellist
\begin{center}
\includegraphics{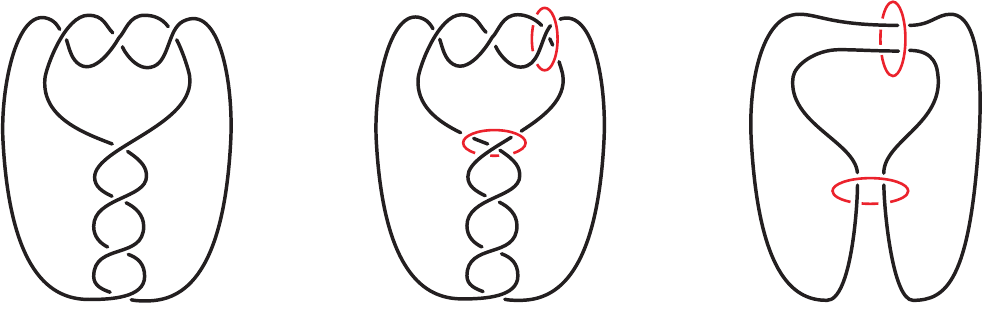}
\caption{The construction of an augmented $2$--bridge link $L$. When $ 
\tw(D)=2$, $L$ is the Borromean rings.}
\llabel{fig:augmentation}
\end{center}
\end{figure}

{\bf Upper bound}\qua The proof of the upper bound uses the  
same surgery argument that Lackenby, Agol, and Thurston used for  
general alternating links \cite{lack-volume}, and the improved  
estimate comes from the special structure of $2$--bridge links.

Let $D$ be a diagram as in \fullref{fig:pillow-braid}. Recall that  
each syllable of $\Omega$ corresponds to a \emph {twist region} where  
two strands of the braid wrap around each other. For every twist  
region, we add an extra link component (called a \emph{crossing  
circle}) encircling the two strands of $K$, obtaining a hyperbolic  
link $J$ \cite{lack-volume}. (See \fullref{fig:augmentation}.)  
Every crossing circle of $J$ bounds a \emph{crossing disk\/} that is  
punctured by the two strands of $K$. Because twice-punctured disks  
are totally geodesic \cite{adams:3ps}, we can untwist all the  
crossings in the twist region and obtain a new link $L$, called an  
\emph{augmented link\/}, whose volume is equal to that of $J$.

When $K$ is a two-bridge link, $L$ has the following alternate  
description. Start with $(\tw(D)-1)$ copies of the Borromean rings, cut  
each one along a crossing disk, and glue the copies together in a  
linear fashion. Volume is additive under this  
operation \cite{adams:3ps}. Thus $\mathrm{Vol}(\mathbb{S}^3 \setminus L) = 2v_8 (\tw(D)-1)$,  
since the Borromean rings have volume $2v_8$. Since $K$ is obtained  
by Dehn filling the crossing circles of $J$, we have
$$\mathrm{Vol}(\mathbb{S}^3 \setminus K) \; < \;  \mathrm{Vol}(\mathbb{S}^3 \setminus  
J) \; =\;
\mathrm{Vol}(\mathbb{S}^3 \setminus L) \; =\;  2v_8 (\tw(D)-1).$$
By choosing a link with many crossings in each twist region, one can  
get $\mathrm{Vol}(\mathbb{S}^3 \setminus K)$ arbitrarily close to this upper  
bound.
\end{proof}

\subsection{Hyperbolic geodesics seen in the projection plane}

For any link diagram $D$, a \emph{crossing arc\/} is a segment  
perpendicular to the projection plane that connects the upper strand  
of a crossing to the lower strand of the same crossing. Morwen  
Thistlethwaite has conjectured that in any reduced alternating  
diagram of a hyperbolic link $K$, every crossing arc is isotopic to a  
hyperbolic geodesic. As a consequence of \fullref{thm:hyperbolic-links}, we can prove this in the case of two-bridge links.

\begin{theorem}\label{thm:crossing-arcs}
Let $D$ be a reduced alternating diagram of a hyperbolic two-bridge  
link $K$. Then every crossing arc of $D$ is isotopic to an edge in  
the Sakuma--Weeks triangulation of $\mathbb{S}^3 \setminus K$, and thus to a  
geodesic.
\end{theorem}

In fact, still more is true: each edge of  
the triangulation is dual to a face of the Ford--Voronoi domain of  
$\mathbb{S}^3 \setminus K$ \cite{aswy-announce}.

\begin{proof}
We begin by observing that the statement is true for a diagram $D_0$  
as in \fullref{fig:pillow-braid}, constructed from an alternating  
braid with one free strand. Every crossing of $D_0$ corresponds to a  
$4$--punctured sphere pleated along edges of our triangulation. As  
\fullref{fig:near-crossing} illustrates, the crossing arc of any  
crossing is isotopic to one of the edges.

As it turns out, the diagram $D_0$ is not overly special.  
Thistlethwaite has proved that every reduced alternating diagram $D$  
of a two-bridge link is standard: that is, $D$ can also be  
constructed from an alternating $4$--string braid, although not  
necessarily with a free strand \cite[Theorem 4.1]{thistle:algebraic}.  
Furthermore, by the Menasco--Thistlethwaite flyping theorem \cite 
{mt:flyping}, we can get from $D_0$ to $D$ by performing a sequence  
of \emph{flypes\/} along the pillowcases of the braid. (See \fullref 
{fig:flype}.) During each flype, the diagram loses a crossing whose  
crossing arc is isotopic to an edge $e$ of the triangulation, and  
gains another crossing, whose arc is isotopic to an edge $e'$ of the  
same slope as $e$. Thus the crossing arcs of every diagram are  
isotopic to geodesics.
\end{proof}

\def\leftmark{F Gu\'eritaud, D Futer}
\def\rightmark{Canonical triangulations}
\bibliographystyle{gtart}
\bibliography{link}

\end{document}